\documentclass[sn-mathphys]{sn-jnl}

\usepackage[utf8]{inputenc}
\usepackage{graphicx,overpic,amsmath,amsfonts,amssymb}
\usepackage{bm}
\usepackage{xcolor}

\def\ba{\begin{aligned}}

\def\be{\begin{equation}}

\def\ea{\end{aligned}}
\def\e{{\bm e}}
\def\ee{\end{equation}}

\newcommand{\anna}[1]{{\color{cyan}#1}}

\title[Surrogate Modeling for NSE with DMD and Manifold Interpolation]{A Data-Driven Surrogate Modeling Approach for Time-Dependent Incompressible Navier-Stokes Equations with Dynamic Mode Decomposition and Manifold Interpolation}

\author*[1]{\fnm{Martin W.} \sur{Hess}}\email{mhess@sissa.it}

\author*[2]{\fnm{Annalisa} \sur{Quaini}}\email{quaini@math.uh.edu}

\author*[1]{\fnm{Gianluigi} \sur{Rozza}}\email{gianluigi.rozza@sissa.it}

\affil[1]{\orgdiv{mathLab}, \orgname{SISSA}, \orgaddress{\street{via Bonomea}, \city{Trieste}, \postcode{34136}, \country{Italy}}}

\affil[2]{\orgdiv{Department of Mathematics}, \orgname{University of Houston}, \orgaddress{\city{Houston, Texas}, \postcode{77204}, \country{USA}}}

\begin{document}

\abstract{
This work introduces a novel approach for data-driven model reduction of time-dependent parametric partial differential equations. Using a multi-step procedure consisting of proper orthogonal decomposition, dynamic mode decomposition and manifold interpolation, the proposed approach allows to accurately recover field solutions from a few large-scale simulations.
Numerical experiments for the Rayleigh-B\'{e}nard cavity problem show the effectiveness of such multi-step procedure in two parametric regimes, i.e.~medium and high Grashof number.
The latter regime is particularly challenging as it nears the onset of turbulent and chaotic behaviour. 
A major advantage of the proposed method in the context of time-periodic solutions is the ability to recover frequencies that are not present in the sampled data.
}



\maketitle

\section{Introduction}

Surrogate modeling, also known as reduced order modeling (ROM), is an invaluable tool for parameter studies of complex dynamical systems that has gained widespread use in recent decades (see \cite{HBMOR_vol1, HBMOR_vol2, HBMOR_vol3}). In this work, we use a non-intrusive (i.e., data-driven) ROM approach, in the sense that only the field solutions of the equations governing the dynamical system at different time steps and parameter values are used to compute the surrogate model. 
Proper orthogonal decomposition (POD), dynamic mode decomposition (DMD), and manifold interpolation are combined into a novel multi-step approach, which allows to recover field solutions at parameters of interest.
As is common for ROM methods, our approach adopts the offline-online decomposition. 
This means that during a time-intensive offline phase all quantities needed for a fast evaluation of solutions over the parameter range
are pre-computed from a few high-fidelity sample solutions. The offline phase can be performed on a high performance cluster, for example. The online phase, which computes the solution for parameters of interest that are not among the sample parameters, can be performed on a laptop or tablet. 

To test and validate our ROM approach, we choose
the Rayleigh-B\'{e}nard cavity problem with fixed aspect ratio and variable Grashof number (Gr), i.e.,~the nondimensional number that describes the ratio of buoyancy forces to viscous forces.  Although this problem features only one physical parameter, it exhibits a wide range of patterns. 
At low Grashof numbers, the system has
unique steady-state solutions. As Gr is increased, the system undergoes several Hopf bifurcations and multiple solutions arise for the same value of the Grashof number. Such solutions past the Hopf bifurcations are time-dependent: they are time-periodic at medium Grashof numbers and 
exhibit turbulent, chaotic behaviour at very high Gr.
A particular difficulty in applying a ROM approach to the Rayleigh-B\'{e}nard cavity over a large range of Gr is the following: the frequencies of time-periodic solutions at online parameters of interest are different from the frequencies at the sample solutions. We have tried several existing ROMs to address this difficulty and have not been successful. This motivated the work presented in this paper. The particular methodology we propose is targeted to problems featuring one or more Hopf bifurcations in the parameter domain of interest, the Rayleigh-B\'{e}nard cavity flow being one challenging example of such problems.

In the setting of bifurcating solutions, ROMs were first considered in \cite{NOOR1982955,Noor:1994,Noor:1983,NOOR198367} for buckling bifurcations in solid mechanics. More recently, in \cite{Maday:RB2,PLA2015162} a reduced basis method is used to track solution branches from bifurcation points arising in natural convection problems. Reduced basis methods are also used in \cite{PR15} to investigate Hopf bifurcations in natural convection problems and in \cite{PITTON2017534} for symmetry breaking bifurcations in contraction-expansion channels. 
Recent works consider ROMs for bifurcating solutions in structural mechanics \cite{PichiRozza2019,pichi2020optcntrl,pichi2021fsi} and physics of condensates \cite{PichiQuainiRozza}.
Finally, we would like to mention that machine learning techniques based on sparse optimization have been applied to detect bifurcating
branches of solutions in \cite{BTBK14,KGBNB17} for a two-dimensional laterally heated cavity and Ginzburg-Landau model, respectively.

The work in this paper builds on our prior work \cite{HessQuainiRozza2022_ETNA,Hess2019CMAME} and
focuses on time-dependent solutions at higher Grashof number than previously investigated. 
At first, we tried the same approach as in 
\cite{HessQuainiRozza2022_ETNA}, which uses artificial neural networks (ANNs) with multilayer perceptrons and different activation functions (like, e.g., ReLU) to improve the localized ROMs introduced in \cite{Hess2019CMAME}. See also \cite{POD_NN_Hesthaven_Ubbiali,pichi2021artificial} for more on POD with ANNs.
However, it turned out that the time evolution of POD coefficients could not be well represented by
 this widely used class of ANNs.
Then, we tried neural ODEs \cite{10.5555/3327757.3327764} and sparse identification of nonlinear dynamics \cite{Brunton3932}, but still failed to recover the correct dynamics.

\begin{figure}[ht!]
\begin{center}
 \includegraphics[scale=.4]{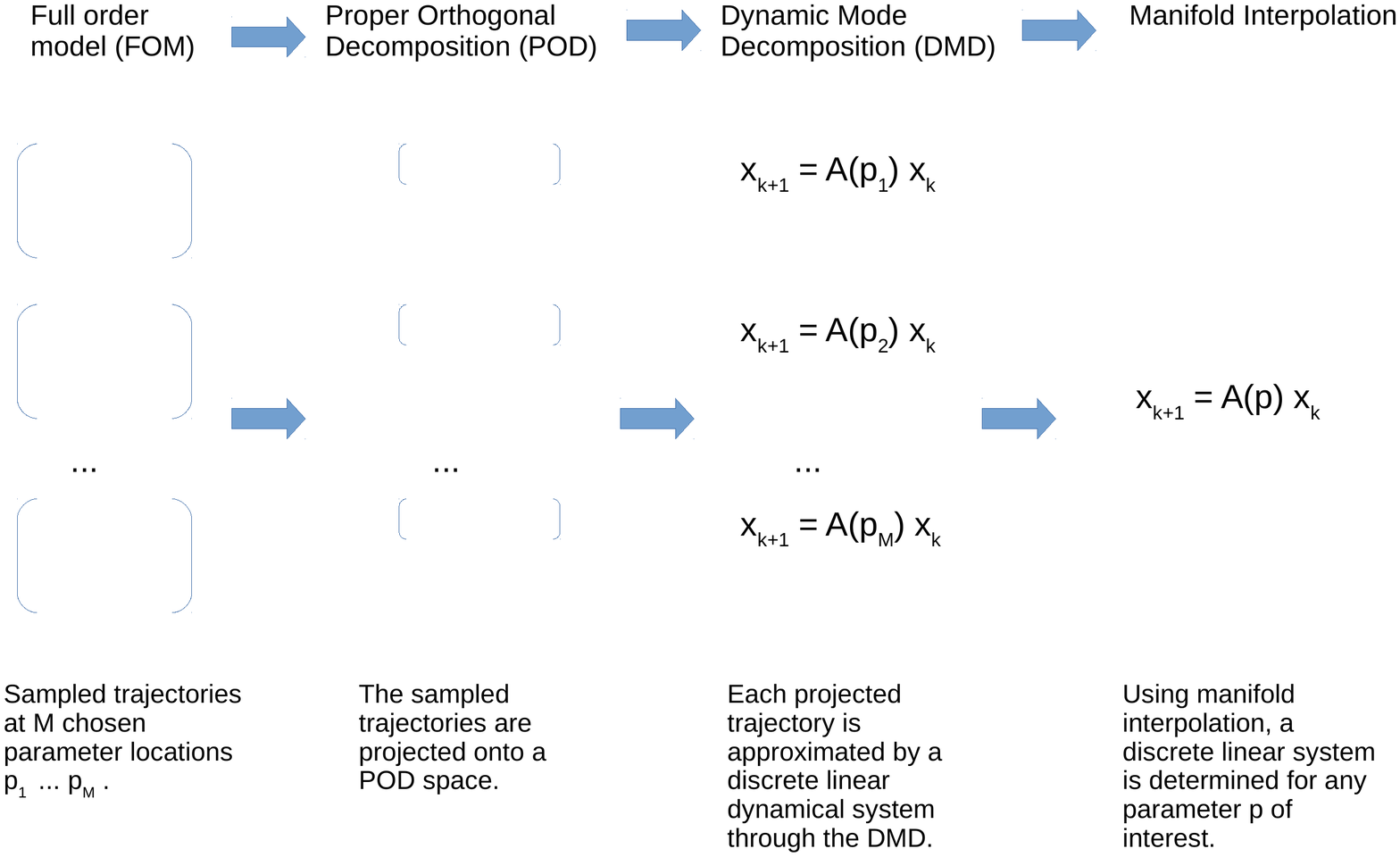} 
 \caption{Schematic of the proposed three-stage ROM. 
 }
 \label{Hess:sketch}
\end{center}
\end{figure}

A major obstacle during the online phase is the correct interpolation of periodicity lengths at intermediate Gr. With increasing Grashof number, the periodicity length of the POD coefficients becomes smaller and the amplitude becomes larger.
The associated flow becomes more complex in general, until it reaches chaotic and turbulent behaviour at very large Gr. In principle, the DMD algorithm (see, e.g., \cite{Koopman1931,doi:10.1137/1.9781611974508,schmid_2010}) is well-suited to resolve time-periodic evolution in a data-driven fashion \cite{schmid_2010}. 
However, the DMD solutions can not be interpolated to intermediate parameter configurations in a straightforward manner. In \cite{GAO2021110907} the DMD is combined with a k-nearest neighbour regression to solve for new parameters of interest, while \cite{doi:10.1063/1.4913868} considers several instances of DMD to solve for parametric problems. Other approaches \cite{TezzeleDemoStabileMolaRozza2020}, \cite{andreuzzi2021dynamic} use different approximation techniques (e.g., polynomial interpolation) and active subspaces to interpolate to new parameters. As mentioned earlier, the issue is that different frequencies are present at intermediate parameters than at the training samples. 
We propose to fix this issue by using 
manifold interpolation based on the DMD operators and DMD modes 
for interpolation at the new parameter values. 
With tangential interpolation based on \cite{Zimmermann2019ManifoldIA}, it is indeed possible to find intermediate frequencies reliably over a wide range of Grashof numbers, which is crucial to accurately recover the time-periodic solutions. 
A schematic of the method is reported in Fig.~\ref{Hess:sketch}.

The rest of the paper is structured as follows. Section 2 introduces the variational formulation of the Navier-Stokes equations governing the Rayleigh-B\'{e}nard cavity and its
discretization with the Spectral Element Method. 
Section 3 explains the model reduction approach and section 4 provides numerical results. In Section 5, 
we provide concluding remarks and further perspectives.


\section{Rayleigh-B\'{e}nard cavity flow}\label{sec:cav}

We consider Rayleigh-B\'{e}nard cavity flow, 
introduced in~\cite{Roux:GAMM} and widely studied since then
(see, e.g. \cite{Gelfgat:Ref11,Hess2019CMAME,PR15})
because of its rich bifurcating behavior, which includes
several Hopf-type bifurcations.
This flow is related to the Rayleigh-B\'{e}nard instability that arises in, e.g., semiconductor crystal growth \cite{KAKIMOTO1995191}. Thus, although simplified, the Rayleigh-B\'{e}nard cavity flow is related to a practical engineering problem.

\subsection{Model description}

The computational domain $\Omega$ is a rectangular cavity with height $1$ and length $4$, i.e.~with aspect ratio $A = 4$, filled with an incompressible, viscous fluid.
The bottom left corner of the cavity is chosen as the origin of the coordinate system. 
The vertical walls are
maintained at constant temperatures $T_0$ (left wall) and $T_0 + \Delta T$ (right wall) with  $\Delta T > 0$, whereas the horizontal walls are thermally insulated. 

This system over a time interval of interest $(0, T)$ is governed by the incompressible \emph{Navier-Stokes} equations
\begin{eqnarray}
\frac{\partial \mathbf{u}}{\partial t} + \mathbf{u} \cdot \nabla \mathbf{u} &=& - \nabla p + \nu \Delta \mathbf{u} + (0, g \beta \Delta T (x/A) \e_y)^T\quad \text{ in } \Omega \times (0, T), \label{Hess:NSE0} \\
\nabla \cdot \mathbf{u} &=& 0\quad \text{ in } \Omega \times (0, T),
\label{Hess:NSE1}
\end{eqnarray}
where $\mathbf{u}$ is the vector-valued velocity, $p$ is the scalar-valued
pressure, and $\nu$ is the kinematic viscosity. Moreover, in \eqref{Hess:NSE0} $g$ denotes the magnitude of the gravitational acceleration, $\beta$ is the coefficient of thermal expansion, 
$x$ is the horizontal coordinate, and $\e_y$ is the unit vector directed along the vertical axis. Problem \eqref{Hess:NSE0}-\eqref{Hess:NSE1}
is endowed with boundary and initial conditions 
\begin{eqnarray}
\mathbf{u} &=& \mathbf{0} \quad \text{ on } \partial \Omega \times (0, T), \\
\mathbf{u} &=& \mathbf{u}_0 \quad \text{ in } \Omega \times 0,
\label{Hess:NSE_boundaryCond}
\end{eqnarray}
with $\mathbf{u}_0$ given. 

The Grashof number
\begin{equation}\label{eq:gr}
\text{Gr} = \frac{g \beta \Delta T}{A \nu^2}
\end{equation}
characterizes the flow regime. 
The Grashof number describes the ratio of buoyancy forces to viscous forces. For large Grashof numbers (i.e., $\gg 1$) buoyancy forces are dominant over viscous forces and vice versa.
Note that with \eqref{eq:gr} we can write the last term in eq.~\eqref{Hess:NSE0} as $(0, \text{Gr}\nu^2 x)^T$.
The Prandtl number for this problem is zero and the viscosity $\nu$ is set to one.

As the Grashof number is increased, the sequence of events is as follows \cite{Roux:GAMM,Gelfgat:Ref11}.
For low Grashof number a steady-state solution exists, which is characterized by a single primary circulation, also referred to as
roll or convective cell.
At a first bifurcation point, the steady-state single roll solution turns into a time-periodic solution and also a steady-state two roll solution appears around the same Gr.
At higher Grashof number, the two roll solutions also turn from steady-state to time-periodic and a three roll steady-state solutions appear.
With increasing Gr, this three roll steady-state solution 
will become time-dependent: 
time-periodic at first and then chaotic (i.e., without an obvious periodicity) upon a further increasing of Gr.
The exact values of the Grashof number where the bifurcations occur 
depend on the aspect ratio $A$ and other parameters, such as the Prandtl number.

\subsection{Discretization}\label{sec:discr}

For the numerical solution of the eq.~\eqref{Hess:NSE0}-\eqref{Hess:NSE1}, 
we adopt the PDE solver \texttt{Nektar++}.
It employs a velocity correction scheme, which advances the nonlinear part explicitly in time and the linear part implicitly.
This time-stepping is also known as a splitting scheme or IMEX (IMplicit-EXplicit) scheme \cite{Guermond_Shen_VCS,Karniadakis_Orszag_Israeli_splitting_methods}. 

\begin{figure}[ht!]
\begin{center}
 \includegraphics[scale=.15]{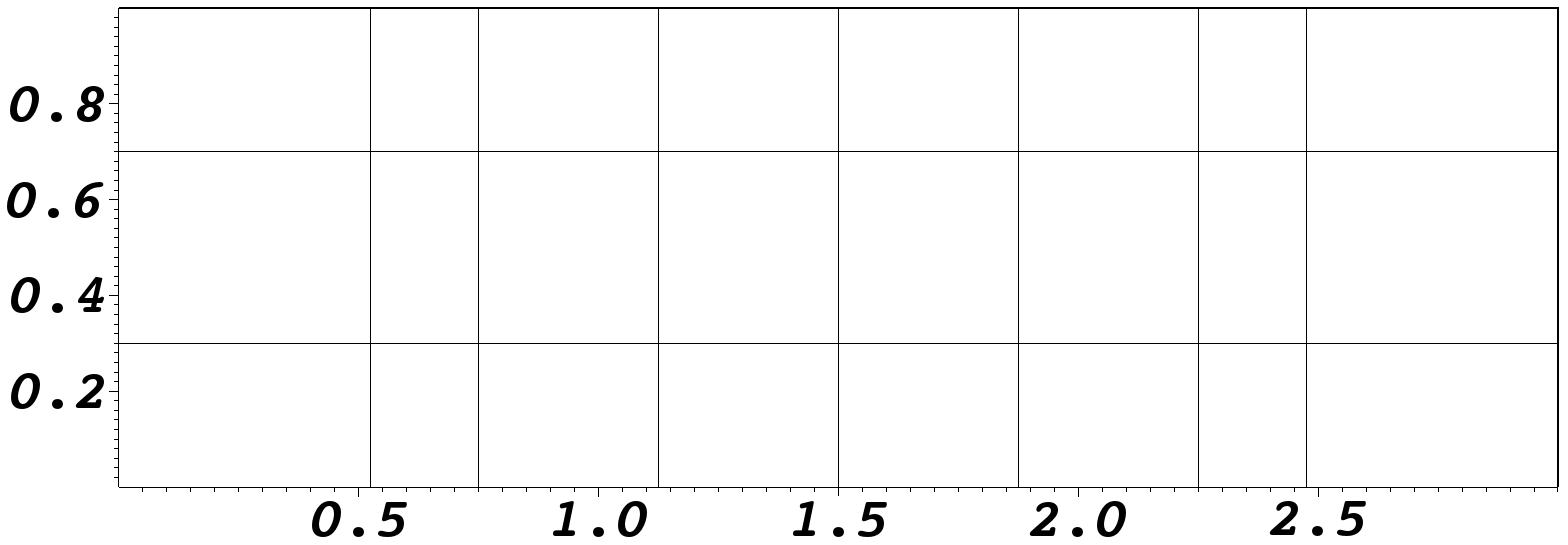} 
 \caption{The $24$ spectral elements used in the simulations, resulting in $18963$ degrees of freedom.}
 \label{Hess:domain_cavity}
\end{center}
\end{figure}

The computational domain is divided into $24$ quadrilateral elements as shown in Fig.~\ref{Hess:domain_cavity}. We use modal Legendre ansatz functions of order $16$, leading to $6321$ global degrees of freedom for each scalar variable (i.e., horizontal velocity, vertical velocity, and pressure ansatz space, this is a standard option in \texttt{Nektar++}), which means a total of $18963$ degrees of freedom.

\begin{figure}[ht!]
\begin{center}
 \includegraphics[scale=.22]{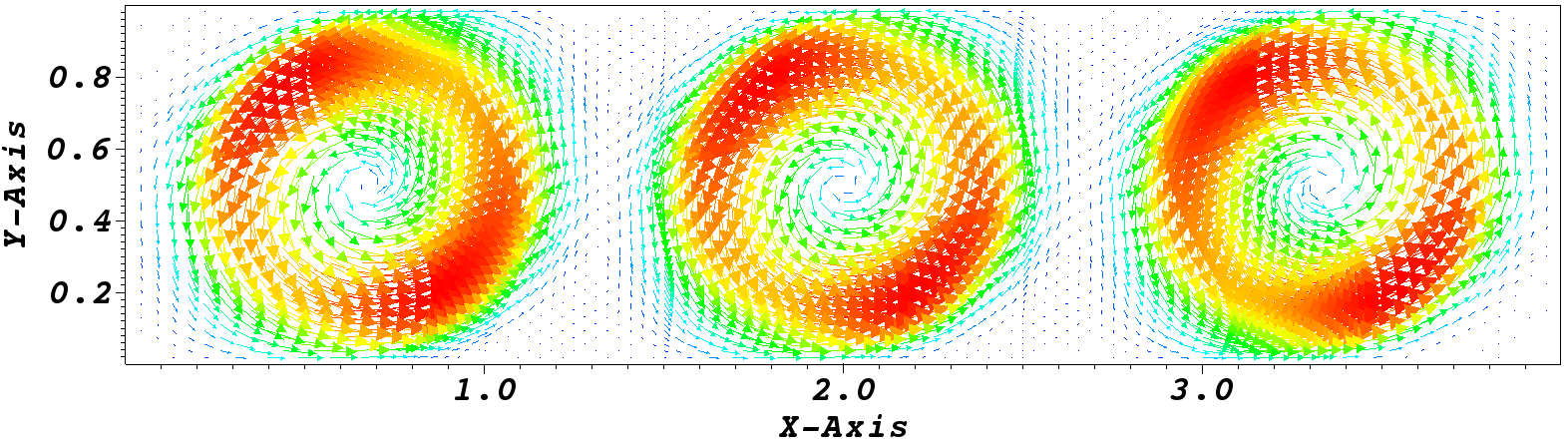} $\quad$
 \includegraphics[scale=.45]{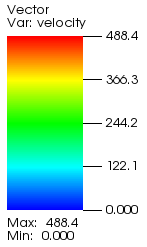}
 \caption{Velocity vector field at Grashof number $100\mathrm{e}{3}$. Color and length of a vector indicate velocity magnitude.}
 \label{Hess:Gr100}
\end{center}
\end{figure}

\begin{figure}
\begin{center}
\begin{overpic}[width=.76\textwidth,grid=false]{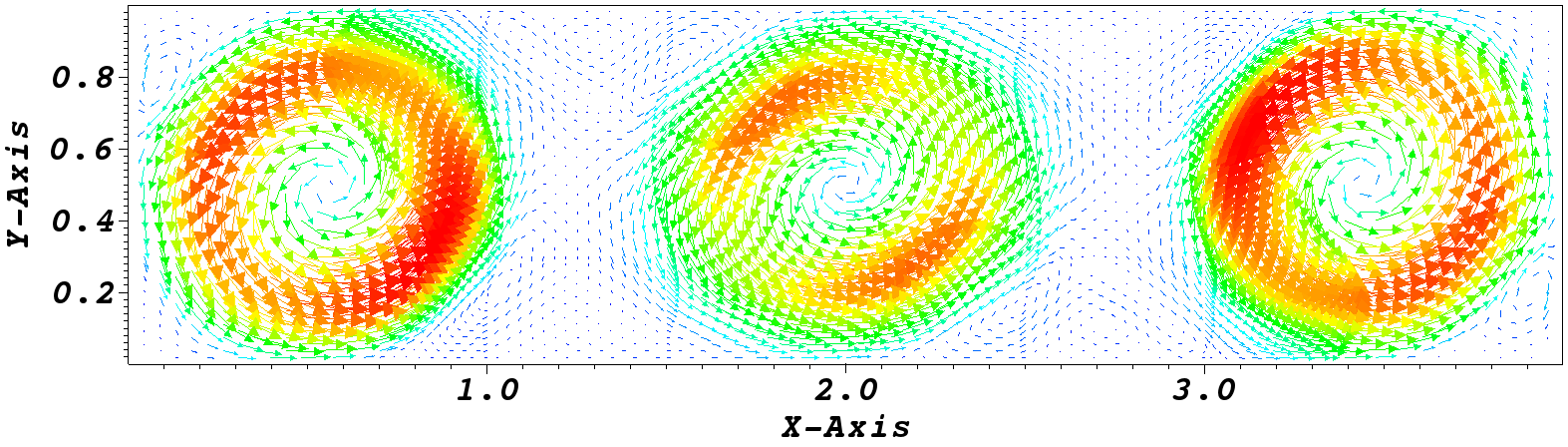}
\put(46,29){\small{$t = 0.705$}}
\end{overpic}$\quad$
 \includegraphics[scale=.45]{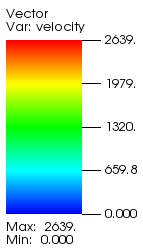} \\
 \begin{overpic}[width=.76\textwidth,grid=false]{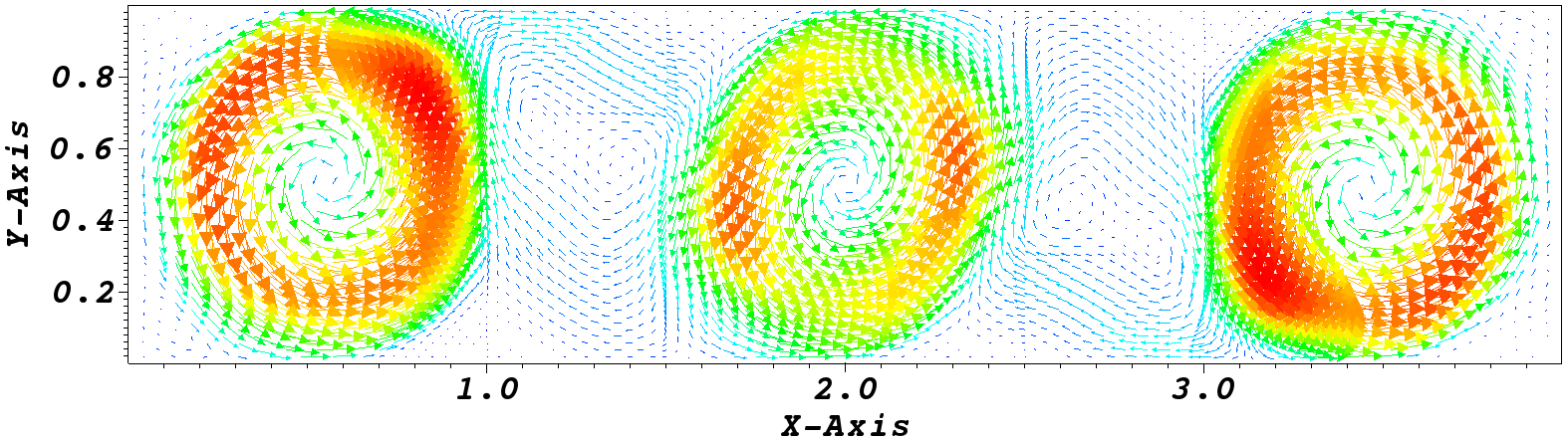}
\put(47,29){\small{$t = 0.71$}}
\end{overpic}$\quad$
 \includegraphics[scale=.45]{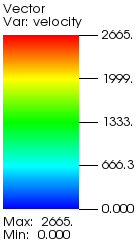} \\
  \begin{overpic}[width=.76\textwidth,grid=false]{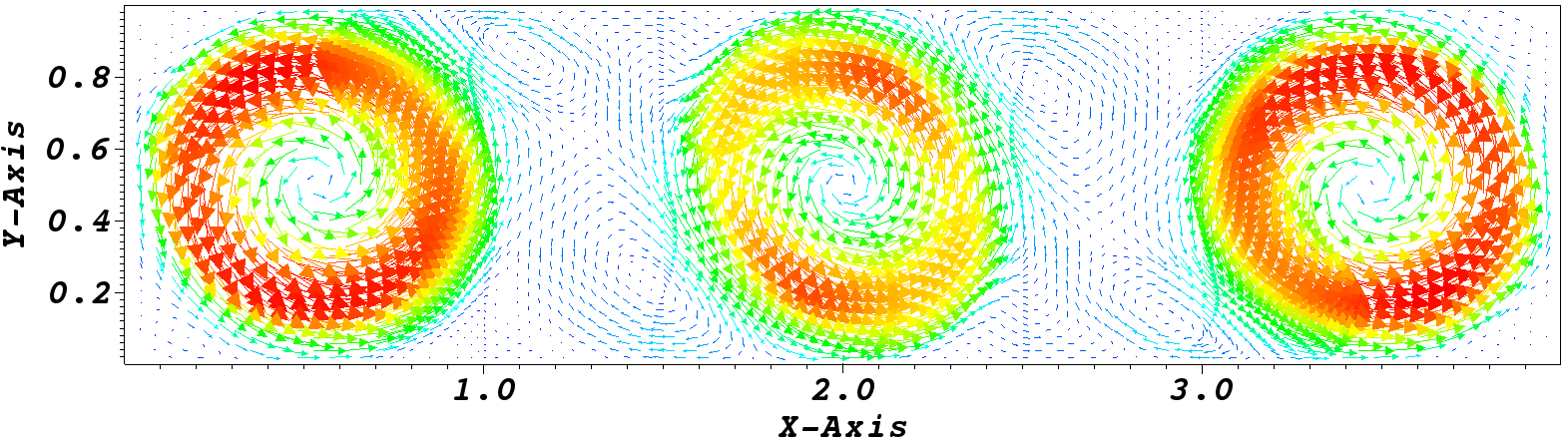}
\put(46,29){\small{$t = 0.715$}}
\end{overpic}$\quad$
 \includegraphics[scale=.45]{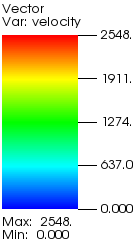} \\
   \begin{overpic}[width=.76\textwidth,grid=false]{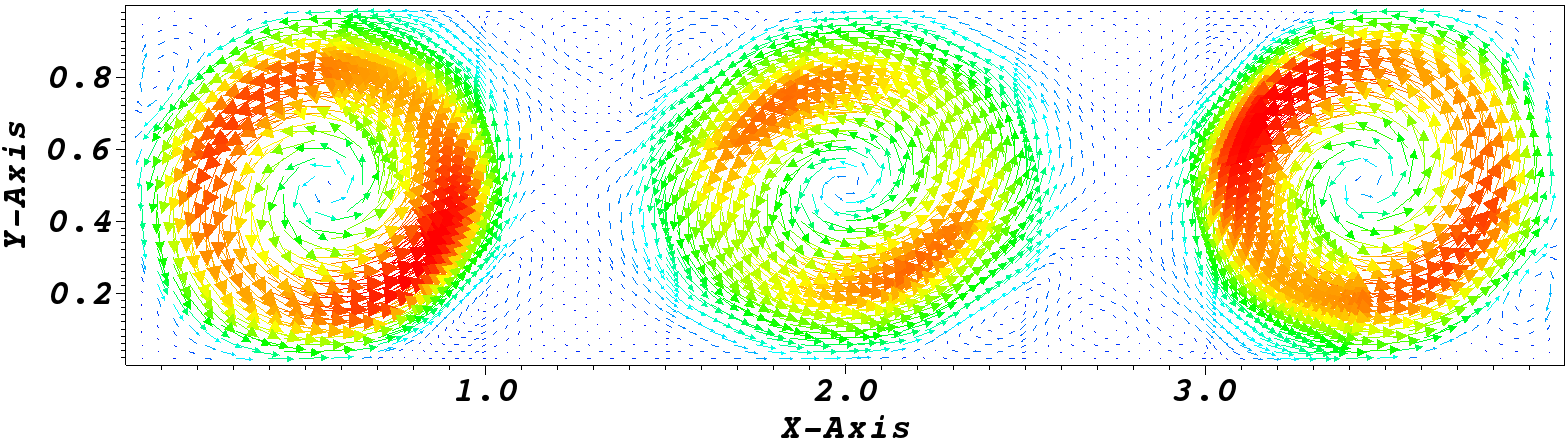}
\put(47,29){\small{$t = 0.72$}}
\end{overpic}$\quad$
 \includegraphics[scale=.45]{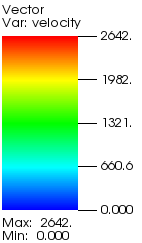} \\
    \begin{overpic}[width=.76\textwidth,grid=false]{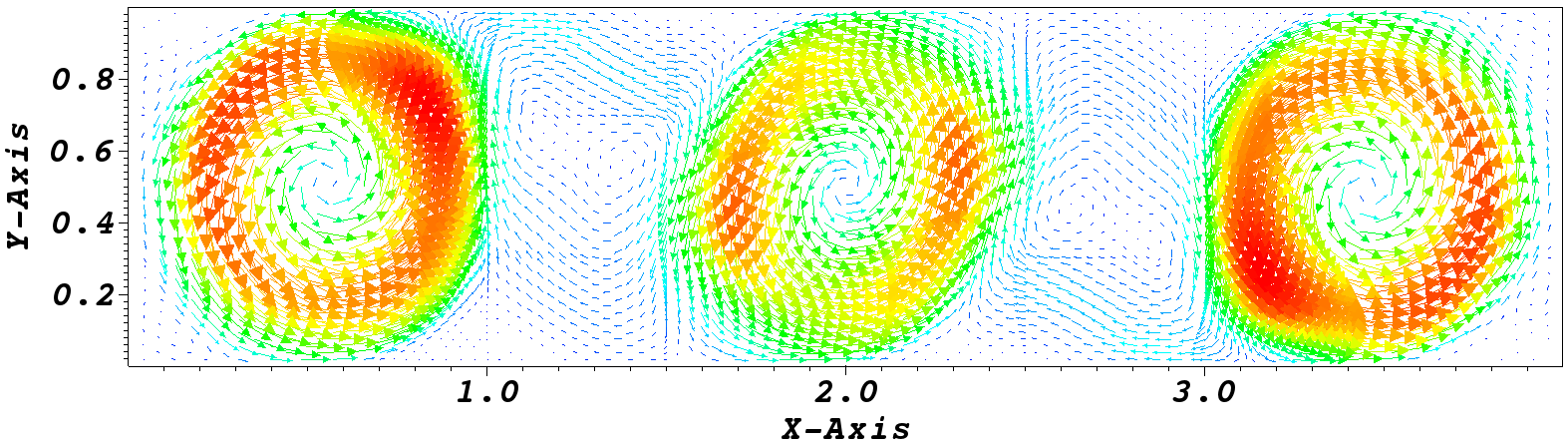}
\put(46,29){\small{$t = 0.725$}}
\end{overpic}$\quad$
 \includegraphics[scale=.45]{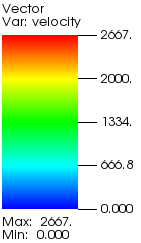} \\
     \begin{overpic}[width=.76\textwidth,grid=false]{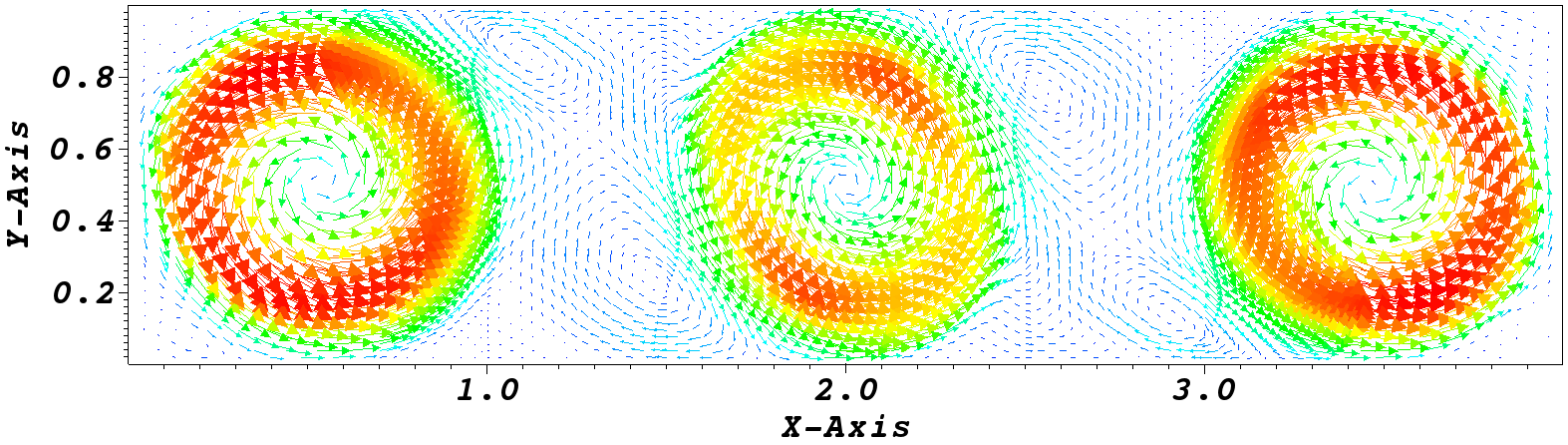}
\put(47,29){\small{$t = 0.73$}}
\end{overpic}$\quad$
 \includegraphics[scale=.45]{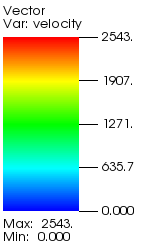} 
 \caption{Time-periodic velocity vector field at Grashof number $700\mathrm{e}{3}$. Color and length of a vector indicate velocity magnitude. 
 Shown are approximately two periods. }
 \label{Hess:Gr700}
\end{center}
\end{figure}

We treat the Grashof number Gr as a parameter and assume it ranges over 
two intervals: $[100\mathrm{e}{3}, 150\mathrm{e}{3}]$ and  $[650\mathrm{e}{3}, 700\mathrm{e}{3}]$.
In both intervals, three roll solutions are typically encountered. 
The velocity vector field at $\text{Gr} = 100\mathrm{e}{3}$ is shown in Fig.~\ref{Hess:Gr100}. 
Upon visual inspection of the flow field, no time-dependence can be detected.
However, it is hard to determine numerically whether this solution is nearing a steady state or is time-periodic since the convergence speed close to the critical value of Gr for the bifurcation point is very slow and a time-periodic pulsation 
with a very small amplitude around a mean field might also be possible.
Fig.~\ref{Hess:Gr700} shows the time-periodic solution at $\text{Gr} = 700\mathrm{e}{3}$ for about two periods. Periodicity is easy to observe from 
the numerical solution. As the Grashof number increases, the period becomes shorter.
For $\text{Gr} > 700\mathrm{e}{3}$, chaotic behaviour can be already observed.
For example, at $\text{Gr} = 1000\mathrm{e}{3}$ we observed that the POD coefficient of the first dominant mode can only be described as "noise", which supports the impression from the velocity time evolution as chaotic. Of course, these are only numerical observations. There is no analytical proof.



Our numerical studies will focus on two distinct parameter domains. First, 
we will look at the interval $[100\mathrm{e}{3}, 150\mathrm{e}{3}]$, where the periods are rather large and the three roll time-dependent solutions have just occurred.
A full-order solution is computed at Gr = $150\mathrm{e}{3}$ over a long time interval to ensure that the limit cycle is reached. Then,
each solution of interest in the interval $[100\mathrm{e}{3}, 150\mathrm{e}{3}]$ is initialized with the solution at Gr = $150\mathrm{e}{3}$. The time
step is set to $1\mathrm{e}{-6}$ and $2\mathrm{e}{5}$ time steps are computed. However, the first $1\mathrm{e}{5}$ time steps are disregarded
in order to  ensure that the solution is sufficiently close to its limit cycle for each parameter of interest. Next, we will consider interval $[650\mathrm{e}{3}, 700\mathrm{e}{3}]$, where the periods are short and the simulations are close to becoming chaotic. Thus, a smaller time step size of $5\mathrm{e}{-7}$ is used. 
We compute $3\mathrm{e}{5}$ time steps and disregard the first $1\mathrm{e}{5}$ time steps. In this second parameter interval, we first compute the full-order solution
at Gr = $650\mathrm{e}{3}$ and use it to initialize all the other solutions of interest.

\section{A model order reduction approach}

The offline phase of our model order reduction approach is articulated into two steps: 
i) proper orthogonal decomposition (POD) briefly explained in Sec.~\ref{sec:POD}
and ii) dynamic mode decomposition (DMD) described in Sec.~\ref{sec:DMD}.
For the online phase, we adopt manifold interpolation as explained in Sec.~\ref{sec:manifold_interp}.

\subsection{Proper Orthogonal Decomposition}\label{sec:POD}

At each Grashof number, we collect 
the velocity field solutions at every time step in the
time interval of interest. These real vectors of dimension $\mathcal{N}$ ($\mathcal{N}$ referring to size of the spatial discretization)
form the trajectory for the given Grashof number. Our first goal is to find a projection matrix to reduce the large dimension $\mathcal{N}$ to a lower dimension $N$.
We achieve this through POD, which computes a projection space $X_N$ used to project the trajectories 
for all parameters in the parameter domain of interest.
Because of the very small time steps required by the cavity simulations, POD is just an initial data reduction step. A second reduction step is needed in order to contain the storage requirements for 
the trajectories. See Sec.~\ref{sec:DMD}.

The POD is based on an operator eigenvalue problem that reduces to the singular value decomposition for discrete data. Given a sample matrix $S \in \mathbb{R}^{\mathcal{N} \times N}$,
compute the singular value decomposition as 
\begin{equation*}
 S = U \Sigma V^T,
\end{equation*}
\noindent where $\Sigma \in \mathbb{R}^{\mathcal{N} \times N}$ is a diagonal matrix with the (non-negative) singular values on the diagonal and $U \in \mathbb{R}^{\mathcal{N} \times \mathcal{N}}$ and $V  \in \mathbb{R}^{N \times N}$ are orthogonal.
Assuming that the singular values are ordered in decreasing order, then the first columns of $U$, called left singular vectors, constitute the dominant POD modes.
The most dominant POD modes are then used as basis functions for the reduced order projection space $X_N$.
For the sake of brevity, we do not report here further details and refer the interested reader to textbooks, such as, \cite{hesthaven2015certified}.

The number of POD modes that are retained is typically determined by a threshold on the percentage of the sum of the singular values, e.g.~$99\%$. 
In particular, if the prescribed threshold is met by the sum of the $L$ largest singular values, but not by the sum of the $L-1$ largest singular values, 
then the $L$ leftmost columns of $U$ are used in the reduced 
order ansatz space $X_N$. See, e.g~\cite{LMQR:2014}, for more details and computational insights on POD in computational fluid dynamics.

The sample matrix $S$ needs to cover the features of the time-dependent solutions over the parameter range
in order for the resulting projection space $X_N$
to retain such features. Because the problem under consideration leads to simulations with large time trajectories, we derive the sample matrix $S$ in an adaptive fashion. For each full-order simulation, we first generate an intermediate matrix by following an adaptive snapshot selection strategy 
from \cite{10.1007/978-3-319-10705-9_42}: we collect samples by adding time instants only if the angle to the already chosen time instants is over a given threshold.
POD is performed on the sampled time instants for each parameter location. Then, the dominant modes resulting from the POD for each parameter location are collected into a second sample matrix, separately for each velocity component. Then, a second application of POD defines the actual projection space $X_N$ and can be understood as a compound POD space of the POD spaces for each time-trajectory. At the end of this first step, we obtain the time-trajectories at each sampled parameter location projected onto the space $X_N$.

This two-tier procedure described in this section allows to keep the storage requirement low, such that the algorithms can even be executed on a common workstation.
In fact, no more than $3$ GB of RAM were necessary to hold a sampled time trajectory.

\begin{figure}[ht!]
\begin{center}
 \includegraphics[scale=.66]{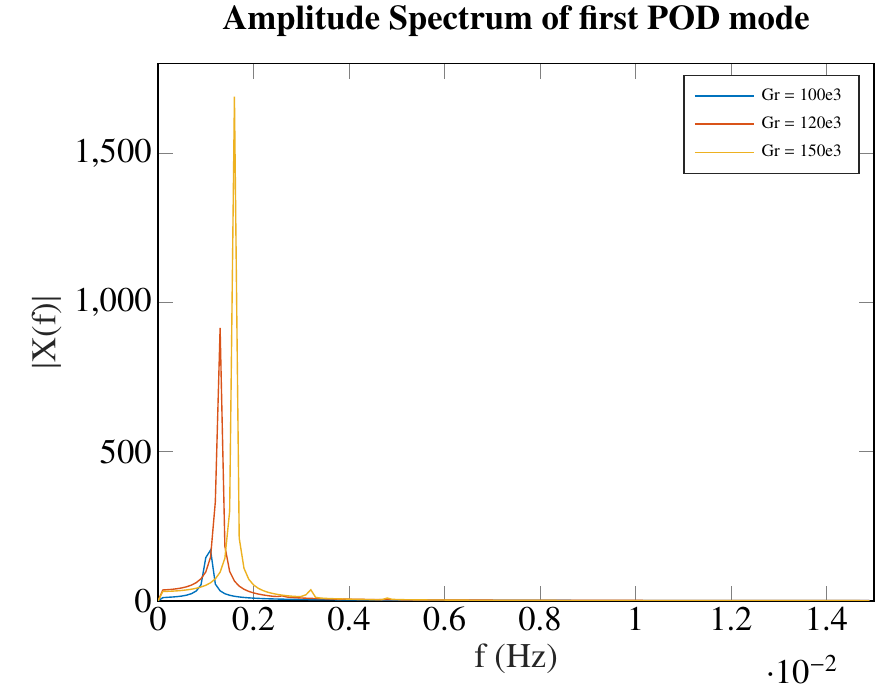} 
 \caption{Amplitude spectrum of the first POD mode for three values of the Grashof number in the medium Grashof interval.}
 \label{Hess:medGr_mode1}
\end{center}
\end{figure}

\begin{figure}[ht!]
\begin{center}
 \includegraphics[scale=.66]{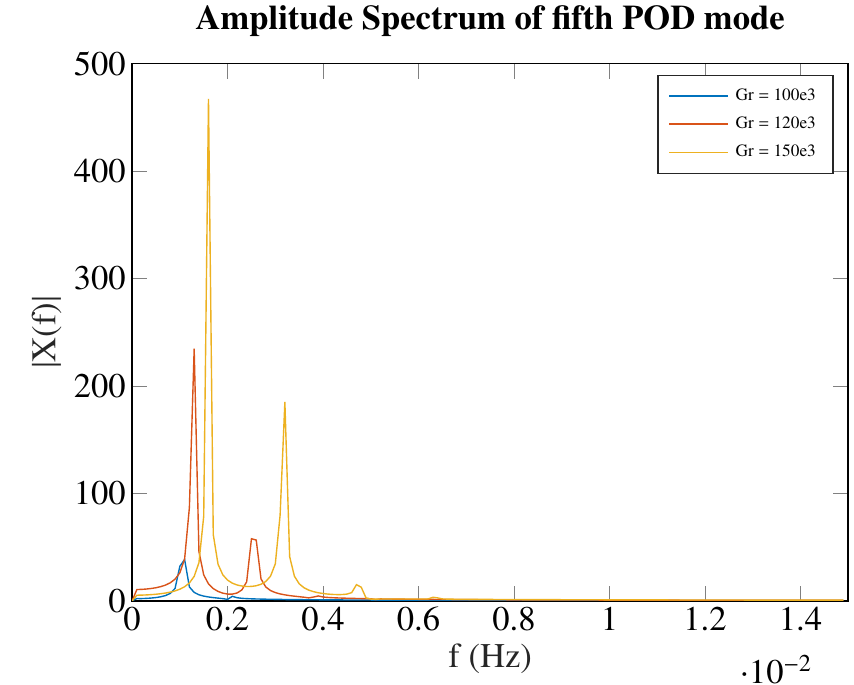} \quad
 \includegraphics[scale=.66]{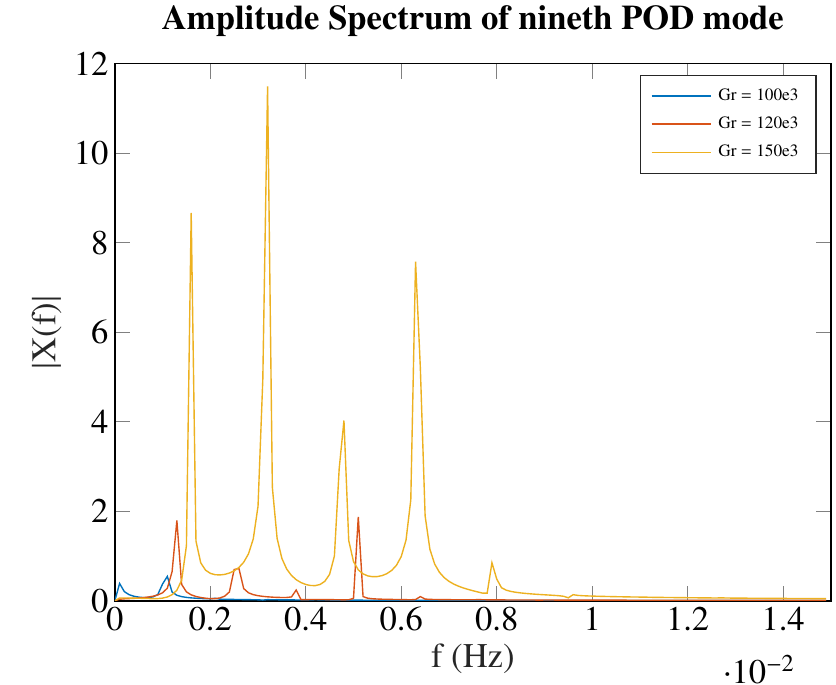} 
 \caption{Amplitude spectrum of the fifth (left) and nineth (right) POD mode for three values of the Grashof number in the medium Grashof interval.}
 \label{Hess:medGr_mode59}
\end{center}
\end{figure}

\begin{figure}[ht!]
\begin{center}
 \includegraphics[scale=.66]{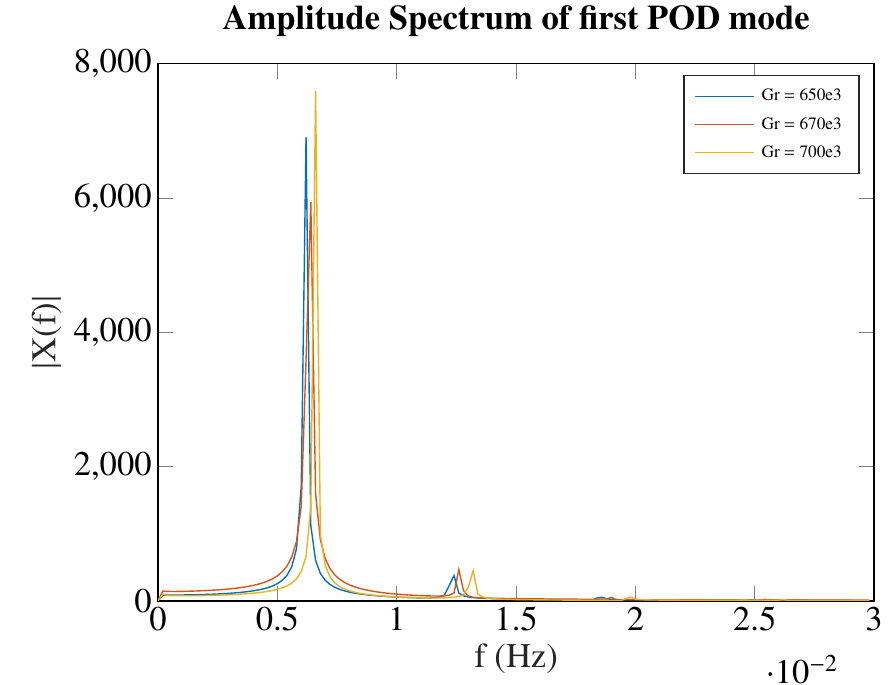} 
 \caption{Amplitude spectrum of the first POD mode for three values of the Grashof number in the high Grashof interval.}
 \label{Hess:highGr_mode1}
\end{center}
\end{figure}

\begin{figure}[ht!]
\begin{center}
 \includegraphics[scale=.64]{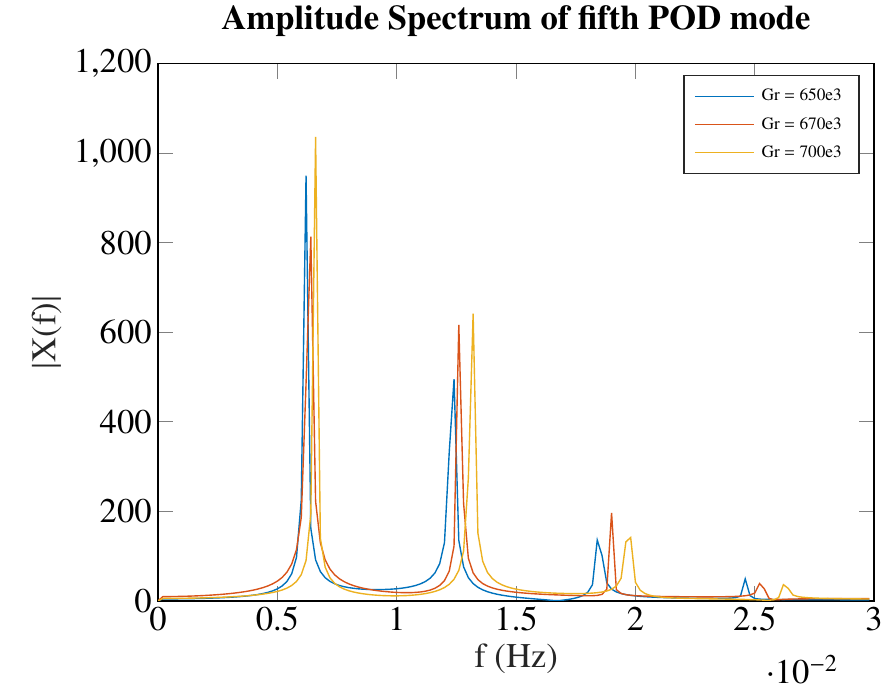} \quad
 \includegraphics[scale=.64]{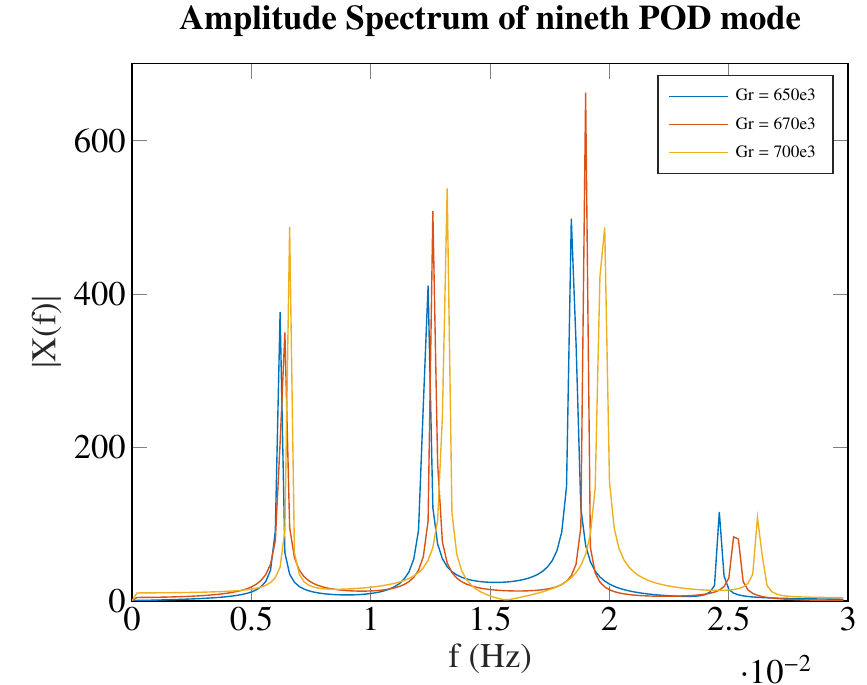} 
 \caption{Amplitude spectrum of the fifth (left) and nineth (right) POD mode for three values of the Grashof number in the high Grashof interval.}
 \label{Hess:highGr_mode59}
\end{center}
\end{figure}

Fig.~\ref{Hess:medGr_mode1} displays the spectrum of the first POD mode of the horizontal velocity component for Grashof numbers $100\mathrm{e}{3}$, $120\mathrm{e}{3}$ and $150\mathrm{e}{3}$.
We observe that the dominant frequency 
increases with increasing Grashof number.
At higher POD modes, more frequencies present, but with a smaller amplitude.
Thus, they are less important for an accurate approximation. See Fig.~\ref{Hess:medGr_mode59} for the amplitude spectrum of the fifth and nineth POD mode of the horizontal velocity component for the same three Grashof numbers.
The same conclusions (i.e., the frequencies increase with increasing Grashof number and more small amplitude frequencies are present in higher POD modes) hold for the high Grashof interval, although in this interval the dominant frequencies are higher and amplitudes are larger than in the medium Grashof interval. See Figs.~\ref{Hess:highGr_mode1} and \ref{Hess:highGr_mode59}.
The amplitude spectra for the POD modes of the vertical velocity component are omitted because they look very similar to Fig.~\ref{Hess:medGr_mode1}-\ref{Hess:highGr_mode59}.

\subsection{Dynamic Mode Decomposition}\label{sec:DMD}

The POD procedure described in the previous section provides a projected trajectory that will take the role that is typically associated with the full-order trajectory in the DMD algorithm. 
We refer the reader to \cite{Koopman1931,doi:10.1137/1.9781611974508,schmid_2010} for an introduction to DMD. For its software implementation, in this work, we use \texttt{PyDMD}\footnote{\textit{https://github.com/mathLab/PyDMD}} \cite{Demo2018}.

Assume the time trajectory is given in the form of state variables $( \mathbf{ x }^k )_{k=1}^m \subset X_N$, 
with $m$ being the total number of time steps.
The goal of DMD is to obtain a linear operator $A \in \mathbb{R}^{N \times N}$, which approximates the dynamics as
\begin{equation}
 \mathbf{ x }^{k+1} \approx A \mathbf{ x }^k \quad \forall k = 1, \ldots, m-1.
 \label{req_DMD}
\end{equation}
If we arrange the state vectors for $k = 1, \ldots, m-1$ column-wise in a matrix X and 
the state vectors for $k = 2, \ldots, m$ column-wise in a matrix Y, then \eqref{req_DMD} is equivalent to
\begin{equation}
 Y \approx A X .
 \label{req_DMD_matrix}
\end{equation}
A best-fit approach computes $A = Y X^\dagger$, where $X^\dagger$ denotes the Moore-Penrose pseudoinverse of $X$.
The linear predictor $A$, also called the \textit{Koopman operator}, allows to recover an approximate trajectory by evaluating \eqref{req_DMD} starting from a given $\mathbf{ x}^1$.

In order to have a reduced order computation of the trajectory, we first compute the rank $r$ truncated singular value decomposition of $X$ as $X \approx U_r \Sigma_r V_r^T$. The matrix $U_r$ holds the real-valued DMD modes as columns.
The reduced operator $A_r \in \mathbb{R}^{r \times r}$ 
is defined as
\begin{equation}
 A_r = U_r^T A U_r = U_r^T Y X^\dagger U_r = U_r^T Y V_r \Sigma_r^{-1} U_r^T U_r =  U_r^T Y V_r \Sigma_r^{-1} ,
 \label{DMD_Ar_composition}
\end{equation}
where we have used the fact that $U_r \in \mathbb{R}^{N \times r}$ is orthogonal. Matrix $A_r$
is used for the reduced order computation of the trajectory as follows:
\begin{equation}
  \mathbf{ x }^{k+1}_r = A_r \mathbf{  x }^k_r .
 \label{ROM_DMD_trajec}
\end{equation}
The full-order trajectory can be approximately recovered as $\mathbf{ x }^k = U_r \mathbf{ x }^k_r$.


There are many variants of DMD for different purposes.  In this work, we use the real-valued standard DMD as shown in eq.~\eqref{DMD_Ar_composition}-\eqref{ROM_DMD_trajec}. In fact, since the initial values of the provided trajectory samples are either on the limit cycle or close to it, the standard DMD is sufficient for an accurate reconstruction of the dynamics. However, if the interest is in recovering the trajectories from a common initial value for all parameters, then a variant of the DMD such as high-order DMD (\cite{LeClainche2017.bib}) 
or Hankel-DMD (\cite{doi:10.1137/17M1125236}) 
can be used. See the \texttt{PyDMD} website for implementations and more details.

\subsection{Manifold interpolation}\label{sec:manifold_interp}

During the online phase, one needs to evaluate the trajectory at a 
new parameter of interest. For this, we have to interpolate 
the reduced DMD operator, which requires 
a structure-preserving interpolation on nonlinear matrix manifolds.
Manifold interpolation has been applied to many problems. See, e.g., \cite{doi:10.1137/100813051,doi:10.1137/130932715,https://doi.org/10.1002/fld.2089,FarhatGrimbergManzoniQuarteroni+2020+181+244,LoiseauBruntonNoack+2020+279+320,PanzerMohringEidLohmann+2010+475+484,doi:10.1137/130942462, GIOVANIS2020113269}. Here, we briefly recapitulate
the basics of manifold interpolation following \cite{Zimmermann2019ManifoldIA}.

As explained in Sec.~\ref{sec:DMD}, the DMD provides a reduced order representation of a trajectory at a fixed Grashof number.
The idea is to sample some trajectories at different Grashof numbers, compute the DMD and then interpolate the Koopman operator $A$ to a new Grashof number of interest.
In particular, the DMD modes $U_r$ and the reduced DMD operator $A_r$ will be interpolated independently\footnote{ Interpolating $A$ directly does not seem promising. However, a possible alternative is to consider the DMD over the complex numbers, if the Riemannian metric is available.} and then matrix $A$ will be obtained using the relation
\begin{equation}
 A = U_r A_r U_r^T .
 \label{Hess:decomposition_Koopman_operator}
\end{equation}

A common picture in reduced basis model reduction is that of the solution manifold, where the solution vectors form a manifold in the ambient discrete space. Similarly, the reduced DMD operators define a manifold  in the space of $r \times r$ matrices. In particular, the reduced DMD operators will be understood as elements the of the general linear group, which forms the manifold $\mathcal{M}$.
Direct interpolation of matrix entries typically lead to poor results. Thus, interpolation is done on the tangent space $T_p^\mathcal{M}$ for a base point $p  \in \mathcal{M}$.
Since the tangent space is flat, a direct interpolation with any interpolation algorithm that expresses the interpolant as a weighted sum of the samples is possible.

Let
\begin{equation*}
 Log_p^\mathcal{M} : \mathcal{M} \mapsto T_p^\mathcal{M} 
\end{equation*}
be the Riemannian logarithm and 
\begin{equation*}
 Exp_p^\mathcal{M} : T_p^\mathcal{M} \mapsto \mathcal{M} 
\end{equation*}
the Riemannian exponential.


For a location $p \in \mathcal{M}$, the interpolation is performed following these steps:
\begin{itemize}
\item[I1.] Given a set of data points $\{p_1, \ldots, p_k\}$, choose first a basis point $p_i$.
\item[I2.] Check that  $Log_{p_i}^\mathcal{M}(p_j)$ is well-defined for all $j = 1, \ldots, k$
and compute $v_j = Log_{p_i}^\mathcal{M}(p_j)$ for all $j$. Here, $v_j = v(\text{Gr}_j)$ where
$\text{Gr}_j$ is the $j^\text{th}$ Grashof number sample location.
\item[I3.] Compute $v^*$ via Euclidean interpolation from the $v_j$, where $v^*$ corresponds to the current Grashof number of interest, and interpolate the matrix entries according to the associated parameters.
\item[I4.] Compute $p^* = Exp_{p_i}^\mathcal{M}(v^*)$ as the interpolated matrix.
\end{itemize}
The above algorithm corresponds to Algorithm 7.1 in \cite{Zimmermann2019ManifoldIA}.

The reduced DMD operator is invertible, so a member of the general linear group of $r \times r$ matrices GL(r).
Since GL(r) is open in the space of all $r \times r$ matrices, the tangent space is simply the space of all $r \times r$ matrices. The simplest choice
for the Riemannian metric is the Euclidean metric, which gives a flat
GL(r). With this choice, the Riemannian exponential of $D$ at a base point $A_r \in GL(r)$ is given by
\begin{equation*}
 Exp_{A_r} (D) = A_r + D
\end{equation*}
and the Riemannian logarithm by
\begin{equation*}
 Log_{A_r} ( \bar{D} ) = \bar{D} - A_r.
\end{equation*}
Other options are possible for the Riemannian metric but will not be contemplated in 
this paper. 

The Grassmann manifold $Gr(N,r)$ is the set of all $r$-dimensional subspaces $\mathcal{U} \subset \mathbb{R}^N$:
\begin{equation*}
 Gr(N,r) = \{ \mathcal{U} \subset \mathbb{R}^N  ~\vert~ \mathcal{U} \text{ subspace},~ \dim(\mathcal{U}) = r \}.
\end{equation*} 
It can be defined as a quotient manifold of the Stiefel manifold 
\begin{equation*}
 St(N,r) = \{ U \in \mathbb{R}^{N \times r} ~\vert~ U^T U = I_r \},
\end{equation*}
through
\begin{equation*}
 Gr(N,r) = St(N,r) / O(r) = \{ [U] ~\vert~ U \in St(N,r) \},
\end{equation*}
where $O(r)$ is the set of the orthogonal $r \times r$ matrices and $I_r$ the $r \times r$ identity matrix. This means that a matrix $U \in St(N,r)$ is the matrix representative of the subspace  $\mathcal{U} \in Gr(N,r)$ if $\mathcal{U} = \text{range}(U)$.
The Grassmann manifold is a typical choice of manifold for projection matrices such as the matrix with POD modes as columns
because the choice of the basis is irrelevant, what matters is the space spanned by the vectors.
Interpolation of the DMD modes is understood as interpolation on the Grassmann manifold.

The composition of the Riemannian exponential and logarithm gives the identity on $Gr(N,r)$. 
However, for the matrix representatives in $St(N,r) / O(r)$ the identity does not necessarily hold. 
See, e.g., \cite{doi:10.1137/100813051} for an explanation on this. 
Thus, 
a modified algorithm for the logarithm is needed for the identity to hold at the matrix level. An example of such
modified algorithm is provided in \cite{Zimmermann2019ManifoldIA}, section 7.4.5.2. It reads as follows: Given a base point representative $U \in St(N,r)$ of $\mathcal{U} = [U] \in Gr(N,r)$ and a point on the manifold
$\widetilde{\mathcal{U}} = [\widetilde{U}] \in Gr(N,r)$ with representative matrix $\widetilde{U} \in St(N,r)$ 
\begin{itemize}
\item[L1.] Compute the SVD of $\widetilde{U}^T U$ as
\begin{eqnarray*}
 \Psi S R^T = \widetilde{U}^T U.
\end{eqnarray*}
\item[L2.] Transition to the Procrustes representative
\begin{equation*}
 \widetilde{U}_* = \widetilde{U} \Psi R^T,
\end{equation*}
and compute the intermediate matrix $L$ as
\begin{equation*}
L = ( I_N - UU^T) \widetilde{U}_*,
\end{equation*}
where $I_N$ is the identity matrix.
\item[L3.] Compute the SVD of $L$
\begin{equation*}
Q \Sigma V^T = L.
\end{equation*}
\item[L4.] Compute the tangent vector 
$\text{Log}^{Gr}_\mathcal{U}(\widetilde{\mathcal{U}})$ on the tangent space $T_\mathcal{U} Gr(N,r)$ as
\begin{equation*}
 \text{Log}^{Gr}_\mathcal{U}(\widetilde{\mathcal{U}}) := Q \arcsin(\Sigma) V^T .
\end{equation*}
\end{itemize}
For a base point representative $U \in St(N,r)$ of $\mathcal{U} = [U] \in Gr(N,r)$ and a tangent vector $\Delta \in T_\mathcal{U} Gr(N,r)$, the exponential computes 
the point $[\widetilde{U}]$ on the Grassmann manifold. 
The algorithm is as follows:
\begin{itemize}
\item[E1.] Compute the SVD of $\text{Log}^{Gr}_\mathcal{U}(\widetilde{\mathcal{U}})$ as
\begin{equation*}
 Q \Sigma V^T = \text{Log}^{Gr}_\mathcal{U}(\widetilde{\mathcal{U}}).
\end{equation*}
\item[E2.] Compute $[\widetilde{U}]$ as
\begin{equation*}
 \widetilde{U} = U V \cos(\Sigma) V^T + Q \sin(\Sigma) V^T.
\end{equation*}
\end{itemize}

\textbf{Remark:}
Since we are dealing with a single parameter, it is possible to use the geodesic interpolation. See  \cite{Zimmermann2019ManifoldIA}, Algorithm 7.2. Geodesic interpolation considers the interval of sampled data points $[\text{Gr}_j, \text{Gr}_{j+1}]$, which includes the unsampled data point $\text{Gr} \in [\text{Gr}_j, \text{Gr}_{j+1}]$. The role of the base point is always taken by the matrix representative, which corresponds to the smaller Grashof number $\text{Gr}_j$.
If the Grashof number of interest is closer to the smaller sampled Grashof number (i.e., $\lvert \text{Gr} - \text{Gr}_j \rvert < \lvert \text{Gr} - \text{Gr}_{j+1} \rvert$), the geodesic interpolation is identical to our approach. However, if the Grashof number of interest is closer to the larger sampled Grashof number (i.e., $\lvert \text{Gr} - \text{Gr}_j \rvert > \lvert \text{Gr} - \text{Gr}_{j+1} \rvert$),  then our approach chooses a different base point.
We observed that the different base point
selected by geodesic interpolation leads to accuracy degrading by a factor of $2-5$.
In \cite{doi:10.2514/1.35374}, the authors did not observe such  sensitivity with respect to the choice of base point. However, there are several differences between our work and \cite{doi:10.2514/1.35374}. One important difference is that we use a global POD with an intermediate DMD, while in \cite{doi:10.2514/1.35374} they interpolate from pre-computed bases at some operating points. The application 
in \cite{doi:10.2514/1.35374}, i.e., aeroelasticity in aircrafts,
is also very different from ours.

\section{Numerical results}

As mentioned in Sec.~\ref{sec:discr}, we will consider
two parameter domains for the Grashof number Gr. The first
domain is $[100\mathrm{e}{3}, 150\mathrm{e}{3}]$ and the associated solutions show the onset of time-dependent three roll flow. 
Since previous works (e.g., \cite{Gelfgat:Ref11,Hess2019CMAME})
deals with lower values of Gr, this first interval is referred to 
as \textit{medium Gr range}. The second domain is $ [650\mathrm{e}{3}, 700\mathrm{e}{3}]$, with associated solutions
that are close to the onset of turbulent and chaotic flow patterns.
We will call this second interval \textit{high Gr range}.
Since the time of a single period decreases with increasing Gr and the flow becomes more complex, we expect that the values of
Gr between the medium and high ranges can also be treated with the presented approach.

\subsection{Medium Grashof range}

The samples in the interval Gr $\in [100\mathrm{e}{3}, 150\mathrm{e}{3}]$ are taken every 
$10\mathrm{e}{3}$, i.e., we collect six samples in total. 
As explained in Sec.~\ref{sec:POD}, we perform 
an adaptive POD for each trajectory and form
the compound POD space for all sample trajectories. 
Both PODs use a threshold of $99.99\%$ of the singular values,
leading to a final dimension of $37$ for the horizontal component of the velocity and $41$ for the vertical component.

The DMD algorithm does not use all the POD modes.
In fact, the DMD uses the first $N = 30$ most dominant POD modes 
in both velocity directions and reduces the dimension to $k=10$ DMD modes. We found that in some cases a further restriction gives more accurate results. 
In particular, we inspected
the first POD mode to check if it shows time-periodic behaviour. This can be seen as an indication of accuracy since we expect to observe time-periodicity. 
If first POD mode is not time-periodic, a second DMD was computed with $N = 10$ modes, which provided accurate results. The two values $N = 30$ and $N = 10$ have been determined empirically.
In total 20 DMDs are computed, i.e., two for each of the ten test samples  for the horizontal and vertical component. The reduction to $N = 10$ modes was applied in three cases.

As mentioned in Sec.~\ref{sec:manifold_interp}, in order
to compute the tangential interpolation we need to choose a base point (step I1). This is a crucial choice since the results are quite sensitive to it. In fact, as one can expect, the interpolation is more accurate the closer the base point is to the online parameter of interest. Thus, for each online test parameter we choose the closest sample point as base point. 
As for the interpolation technique, in principle one can
choose any technique that expresses the interpolant as a weighted sum of the samples. The first obvious choice to try is a linear interpolation between the two closest sample points. We found that linear interpolation gives very accurate results that usually cannot be improved by switching to a higher order interpolation or radial basis function interpolation. Thus, we stuck to linear interpolation.

In order to evaluate the accuracy of the ROM approach,
we select psuedo-randomly ten test points in the interval Gr $\in [100\mathrm{e}{3}, 150\mathrm{e}{3}]$.
By ``pseudo-randomly'' we mean that we ensure that no test point coincides with a sample point and that the test points cover the entire parameter interval.  For each test point, we compute the relative $L^2(\Omega)$ and $L^\infty(\Omega)$ error for the velocity for all time steps with respect to the full order simulation. We run the simulations for a total time $T=1\mathrm{e}{-1}$ with time step $1\mathrm{e}{-6}$. Since the
DMD can not properly resolve the swing-in phase (the first $1\mathrm{e}{5}$ time steps), that is not considered for the 
error computations. 
In this way, a start value close to the limit cycle is provided at each test point.

\begin{table}[ht!]\caption{Mean and maximum relative errors for the velocity over $T=1\mathrm{e}{-1}$. The samples in the interval $[100\mathrm{e}{3}, 150\mathrm{e}{3}]$ are taken every
$10\mathrm{e}{3}$.}
\begin{center}
\begin{tiny}
\begin{tabular}{l*{5}{c}r}
  & Gr$=113.8\mathrm{e}{3}$ & Gr$=133.985\mathrm{e}{3}$ & Gr$=132.755\mathrm{e}{3}$ & Gr$=108.13\mathrm{e}{3}$ & Gr$=105.95\mathrm{e}{3}$  \\
\hline
mean $L^2$  &  0.0181 &    0.0072  &   0.0060  &   0.0070  &   0.0201   \\
mean $L^\infty$ &   0.0222  &    0.0076   &   0.0062   &   0.0088   &   0.0258    \\
max $L^2$  & 0.0399   &     0.0169    &    0.0139    &    0.0249   &     0.0628     \\
max $L^\infty$   &  0.0497  &  0.0193  &    0.0152 &   0.0313  &  0.0798   \\[1ex]
   & Gr$=124.92\mathrm{e}{3}$ & Gr$=147.987\mathrm{e}{3}$ & Gr$=117.02\mathrm{e}{3}$ & Gr$=129.26\mathrm{e}{3}$ & Gr$=111.19\mathrm{e}{3}$  \\
\hline
mean $L^2$  &  0.0153   &   0.0067   &   0.0074   &   0.0039  &    0.0046  \\
mean $L^\infty$ &   0.0192  &     0.0071   &    0.0087    &   0.0041  &     0.0057   \\
max $L^2$ & 0.0359   &     0.0174    &    0.0184   &     0.0104  &      0.0123   \\
max $L^\infty$   &   0.0461 &   0.0207  &  0.0228  &  0.0122 &   0.0153  \\[1ex]
\end{tabular}
\end{tiny}
\end{center}
\label{medium_Gr_results}
\end{table}

\begin{figure}[ht!]
\begin{center}
 \includegraphics[scale=.7]{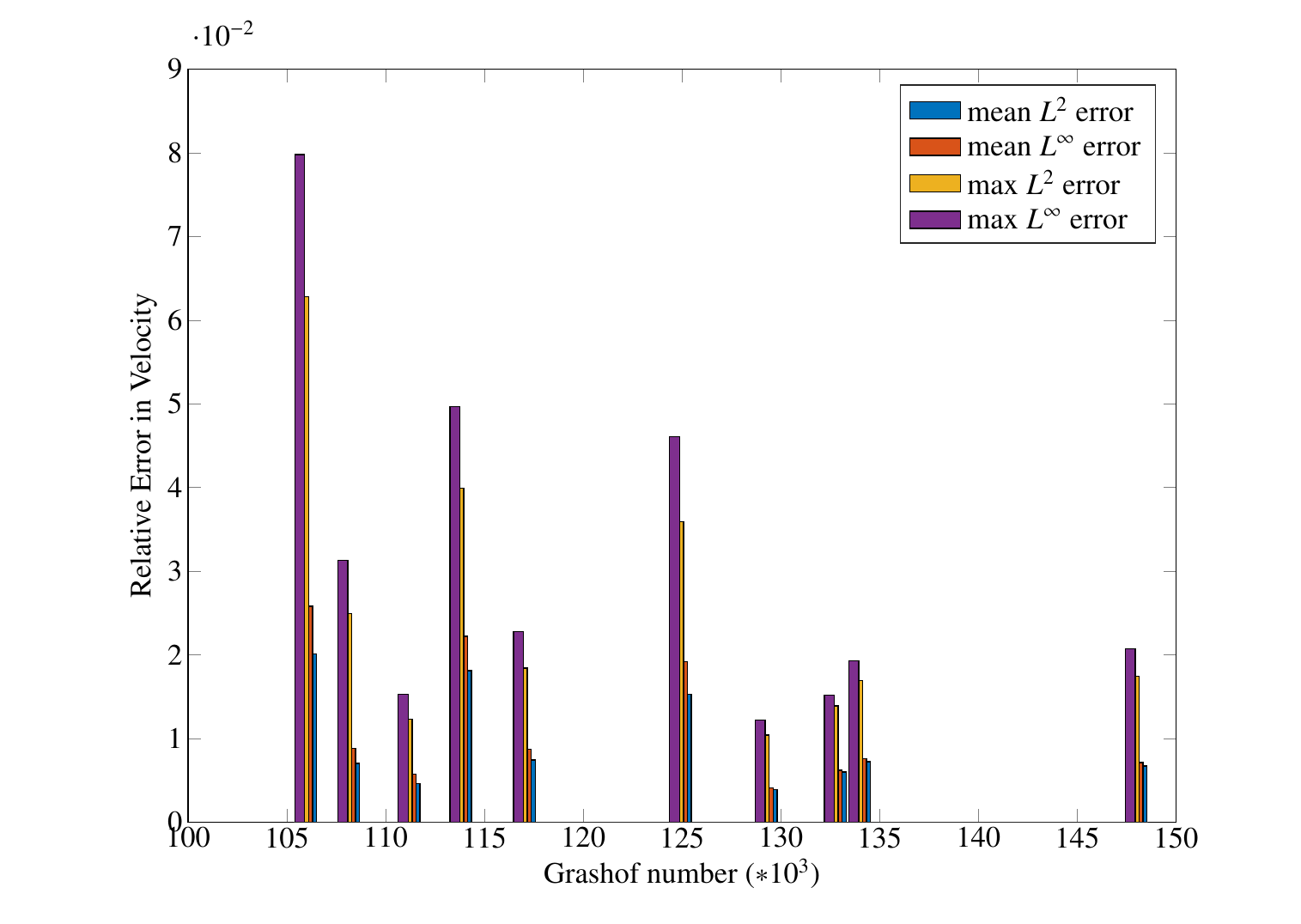} 
 \caption{Bar plot visualizing the results from Table~\ref{medium_Gr_results}.
 The samples in the interval $[100\mathrm{e}{3}, 150\mathrm{e}{3}]$ are taken every
$10\mathrm{e}{3}$.
 }
 \label{medium_Gr_results_bar}
\end{center}
\end{figure}

Table~\ref{medium_Gr_results} reports the mean and maximum 
relative $L^2(\Omega)$ and $L^\infty(\Omega)$ error for the medium Grashof range and Fig.~\ref{medium_Gr_results_bar} visualizes the same data. 
We see that the three test points in-between sample points (Gr$=113.8\mathrm{e}{3}$, Gr$=105.95\mathrm{e}{3}$, and Gr$=124.92\mathrm{e}{3}$) have mean $L^2(\Omega)$ and $L^\infty(\Omega)$ errors up to $2.58\%$. All the remaining test points, which are closer to a sample point, have mean errors below $1\%$. The same observation holds true for the maximum error. In particular, we notice that the maximum $L^\infty(\Omega)$ error for Gr$=105.95\mathrm{e}{3}$ goes up to $8\%$. This shows that the distance to the base point is crucial for the accuracy of our approach, as mentioned above. Thus, we conclude that the proposed interpolation approach provides accurate approximations so long as the sample density is appropriate, i.e.~there is a base point for the manifold interpolation near each new parameter value.

\begin{figure}[ht!]
\begin{center}
 \begin{overpic}[width=0.49\textwidth, grid = false]{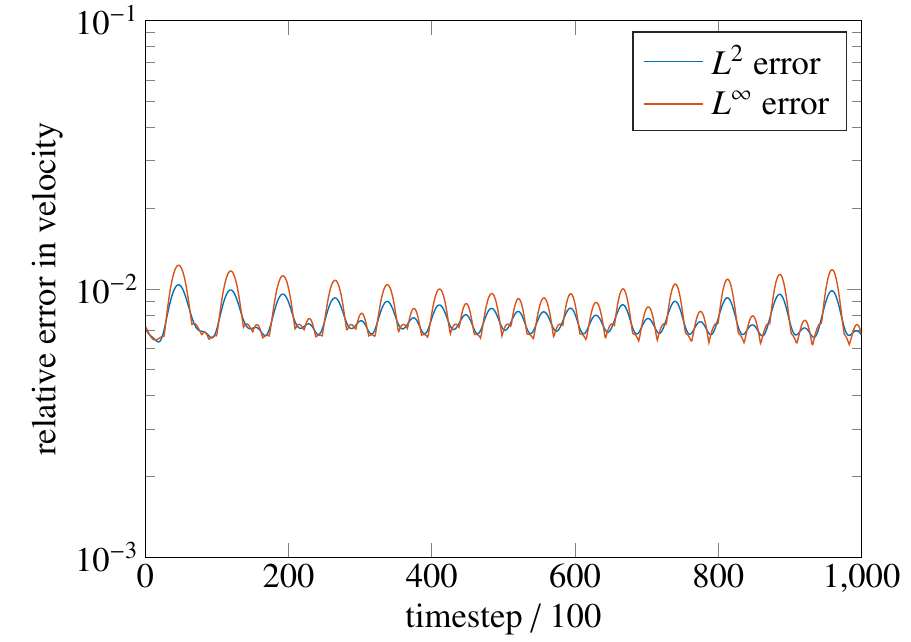}
        \put(27,70){Gr$=129.26\mathrm{e}{3}$ (best)}
      \end{overpic} 
 \begin{overpic}[width=0.49\textwidth, grid = false]{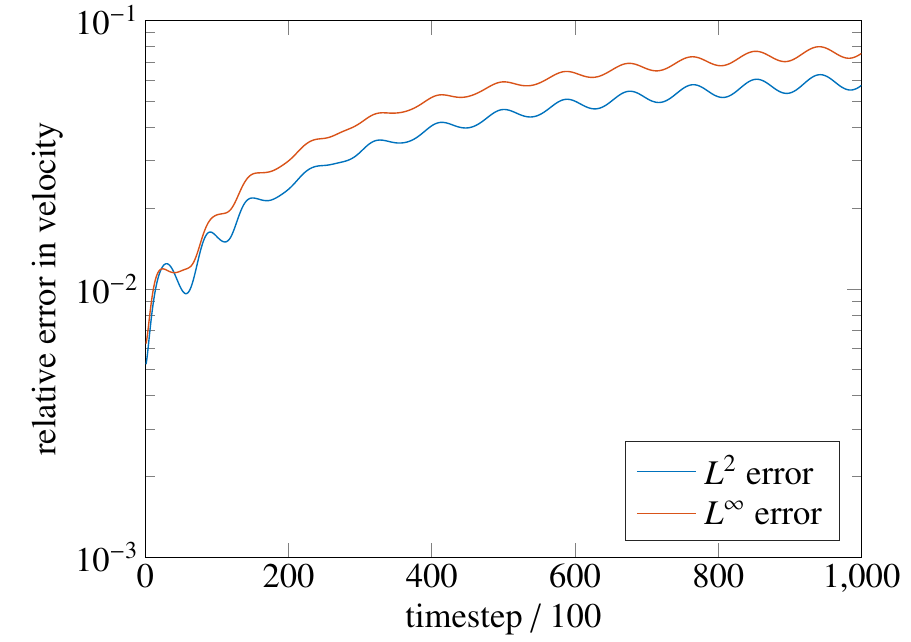} 
        \put(27,70){Gr$=105.95\mathrm{e}{3}$ (worst)}
       \end{overpic} 
 \caption{Relative errors over time for the velocity at Gr$=129.26\mathrm{e}{3}$ (left) and Gr$=105.95\mathrm{e}{3}$ (right). 
 }
 \label{Hess:batch1}
\end{center}
\end{figure}

\begin{figure}[ht!]
\begin{center}
 \begin{overpic}[width=0.49\textwidth, grid = false]{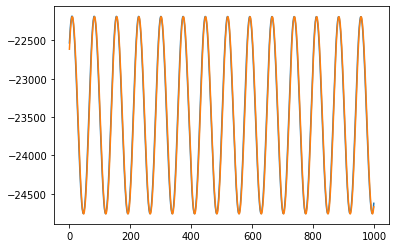}
        \put(30,62){horizontal component}
        \put(-10,35){\small{best}}
      \end{overpic} 
 \begin{overpic}[width=0.49\textwidth, grid = false]{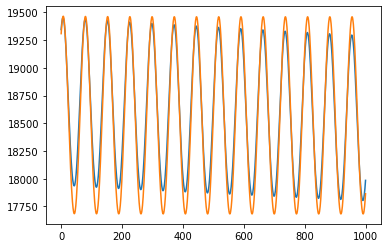} 
        \put(32,63){vertical component}
       \end{overpic}  \\
   \begin{overpic}[width=0.49\textwidth, grid = false]{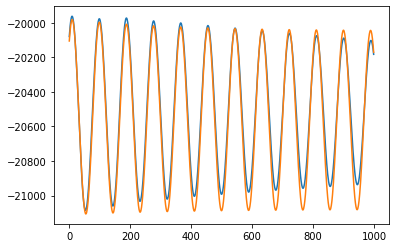}
       \put(-10,35){\small{worst}}
      \end{overpic} 
 \begin{overpic}[width=0.49\textwidth, grid = false]{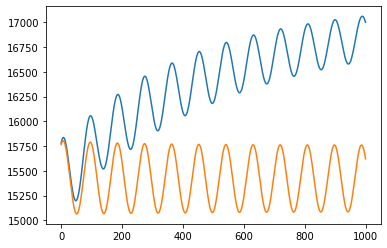} 
       \end{overpic}     
 \caption{Approximation of the first POD mode coefficient (blue) and reference solution (red) for the horizontal (left) and vertical (right) component of the velocity in the best approximation case, i.e.~Gr$=129.26\mathrm{e}{3}$
 (top), and worst approximation case, i.e.~Gr$=105.95\mathrm{e}{3}$ (bottom).}
 \label{Hess:batch1_mode0_xy}
\end{center}
\end{figure}

The relative $L^2(\Omega)$ and $L^\infty(\Omega)$ errors over time for the best approximated case (Gr$=129.26\mathrm{e}{3}$) and the worst approximated case (Gr$=105.95\mathrm{e}{3}$) amongst the test samples in Table~\ref{medium_Gr_results} are shown in Fig.~\ref{Hess:batch1}. 
We see that 
at Gr$=129.26\mathrm{e}{3}$ there is no initial growth of the error over time in contrast to Gr$=105.95\mathrm{e}{3}$.
Let us take a look at the approximation of the first POD modes for the horizontal and vertical component of the velocity for both cases, which are reported in Fig.~\ref{Hess:batch1_mode0_xy}. In the case of Gr$=129.26\mathrm{e}{3}$, the error is dominated by the approximation in the vertical component since the first POD mode is well approximated for the horizontal component. Indeed, the blue line is superimposed to the red line in the top left panel in Fig.~\ref{Hess:batch1_mode0_xy}. Also at 
Gr$=105.95\mathrm{e}{3}$ the first POD mode is
well approximated for the horizontal component, although the difference between approximated
and reference mode becomes more evident as time
passes. In addition, the mismatch between approximated and reference mode for the vertical component of the velocity is much larger at Gr$=105.95\mathrm{e}{3}$ than at Gr$=129.26\mathrm{e}{3}$.
For these two examples, the errors reported in Table~\ref{medium_Gr_results} can be understood  form the approximation of the first POD mode as shown in Fig.~\ref{Hess:batch1_mode0_xy}. For other values of Gr, it is necessary to also look at the other (less dominant) POD modes.

\subsection{High Grashof range}

Following what we have done in the medium Grashof range, we take samples every $10\mathrm{e}{3}$ in the high Grashof range $[650\mathrm{e}{3}, 700\mathrm{e}{3}]$ for a total of six samples.
We repeat the two-tier POD procedure illustrated in Sec.~\ref{sec:POD} and set the threshold for both PODs to $99.99\%$ of the singular values. The final dimensions are $78$ for the horizontal velocity component and $82$ for the vertical velocity component. Notice that the dimensions for both velocity components \anna{are larger} than in the medium Grashof range. 

Just like for the medium Grashof range, the DMD  uses the $N = 30$ of the most dominant POD modes in both velocity directions and reduces to $k=10$ DMD modes. Moreover, if the first POD mode is not showing time-periodic behaviour 
the DMD algorithm is applied again with $N = 10$. 
This was used eight times out of 20 DMDs for the 10 test samples with independent DMDs for the horizontal and vertical component.
Again, the manifold interpolation chooses the closest sample point as base point and uses linear interpolation in the tangent space.

The ten test points are chosen by 
shifting the ten random test points used for medium Gr interval to the high Gr interval Gr $\in [650\mathrm{e}{3}, 700\mathrm{e}{3}]$.
For each test point, we compute the relative $L^2(\Omega)$ and $L^\infty(\Omega)$ error for the velocity for all time steps with respect to the full order simulation.
Once again we remove 
the swing-in phase (the first $1\mathrm{e}{5}$ time steps) from the error computations, so that for each test point a start value close to the limit cycle is provided. The error computation is then performed over another $2\mathrm{e}{5}$ time steps with a time step size of $5\mathrm{e}{-7}$ for a total time $T=1\mathrm{e}{-1}$.

\begin{table}[ht!]\caption{Mean and maximum relative errors for the velocity over $T=1\mathrm{e}{-1}$.
The samples in the interval $[650\mathrm{e}{3}, 700\mathrm{e}{3}]$ are taken every $10\mathrm{e}{3}$.
} 
\begin{center}
\begin{tiny}
\begin{tabular}{l*{5}{c}r}
  & Gr$=663.8\mathrm{e}{3}$ & Gr$=683.985\mathrm{e}{3}$ & Gr$=682.755\mathrm{e}{3}$ & Gr$=658.13\mathrm{e}{3}$ & Gr$=655.95\mathrm{e}{3}$  \\
\hline
mean $L^2$  & 0.1136  &  0.1461 &   0.1197  &  0.0413  &  0.0428     \\
mean $L^\infty$ &  0.1502  &  0.1956  &  0.1592  &  0.0454  &  0.0478   \\
max $L^2$  & 0.1800  &    0.1846   &   0.1421  &    0.0504    &  0.0519      \\
max $L^\infty$   & 0.2428  &   0.2488 &    0.1898 &   0.0575  &   0.0594   \\[1ex]
   & Gr$=674.92\mathrm{e}{3}$ & Gr$=697.987\mathrm{e}{3}$ & Gr$=667.02\mathrm{e}{3}$ & Gr$=679.26\mathrm{e}{3}$ & Gr$=661.19\mathrm{e}{3}$  \\
\hline
mean $L^2$  &     0.0384    &    0.0970    &    0.1039    &    0.0426    &    0.0598 \\
mean $L^\infty$ &    0.0404    &    0.1268    &    0.1365    &    0.0474    &    0.0735 \\
max $L^2$ &  0.0438   &   0.1405   &   0.1681   &   0.0516   &   0.0975  \\
max $L^\infty$   &   0.0472  &     0.1922   &    0.2264   &    0.0596   &    0.1295 \\[1ex]
\end{tabular}
\end{tiny}
\end{center}
\label{high_Gr_results}
\end{table}

\begin{figure}[ht!]
\begin{center}
 \includegraphics[scale=.7]{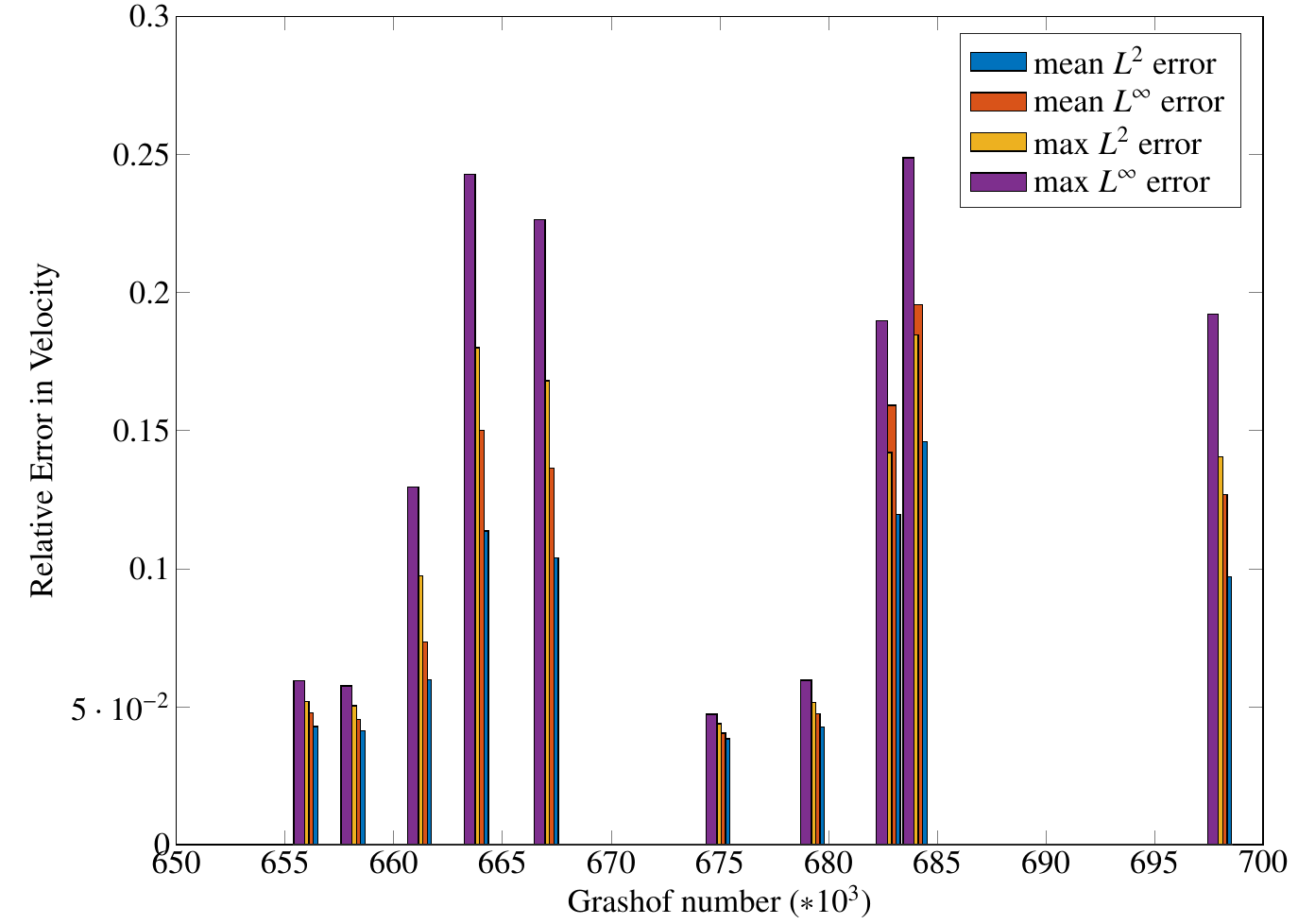} 
 \caption{Bar plot visualizing the results from Table~\ref{high_Gr_results}. The samples in the interval $[650\mathrm{e}{3}, 700\mathrm{e}{3}]$ are taken every $10\mathrm{e}{3}$.}
\end{center}
\label{high_Gr_results_bar}
\end{figure}

Table~\ref{high_Gr_results} reports the mean and maximum relative $L^2(\Omega)$ and $L^\infty(\Omega)$ error for the velocity and Fig.~\ref{high_Gr_results_bar} visualizes the same data. 
From this table, it is not easy to guess when
an approximation is more or less accurate. The distance to the bast point in the manifold interpolation does not seem to play the same obvious role as in the medium Gr range.
For example, the mean $L^2(\Omega)$ error is less than $6\%$ for half the test points, while it goes up to about $15\%$ for Gr$=683.985\mathrm{e}{3}$. Similar observations can be made for the mean $L^\infty(\Omega)$ error and the maximum errors.

\begin{figure}[ht!]
\begin{center}
 \begin{overpic}[width=0.49\textwidth, grid = false]{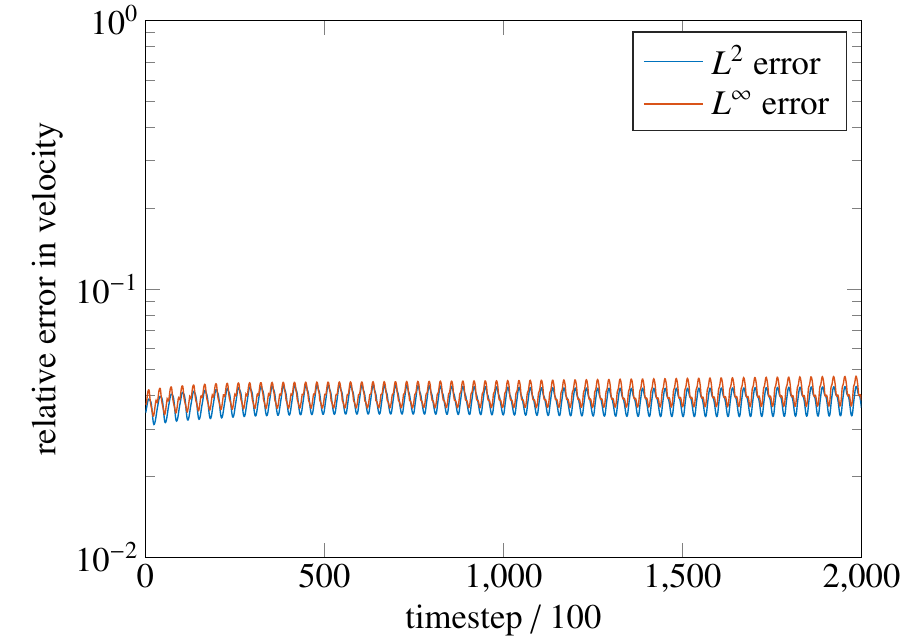}
        \put(27,70){Gr$=674.92\mathrm{e}{3}$ (best)}
      \end{overpic} 
 \begin{overpic}[width=0.49\textwidth, grid = false]{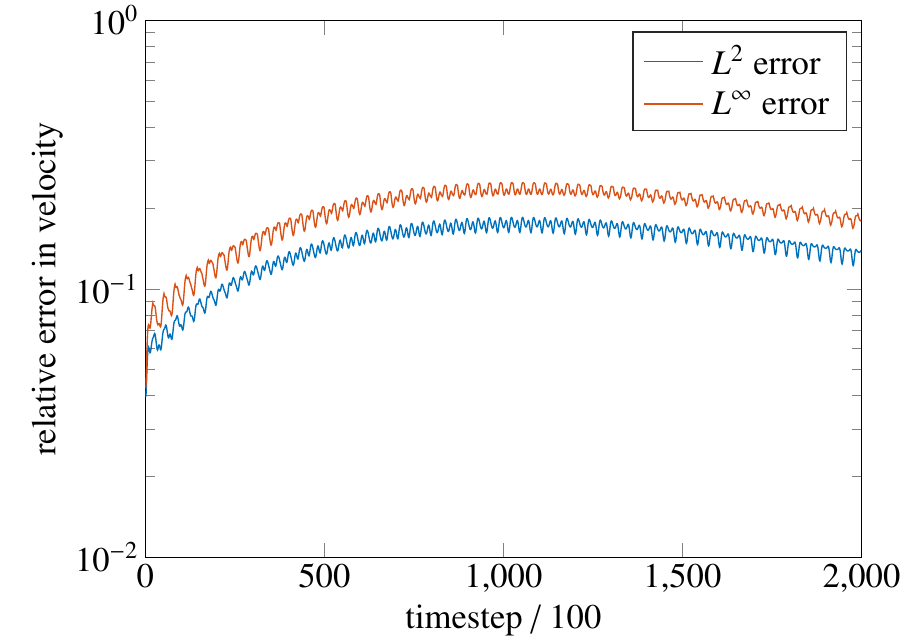} 
        \put(27,70){Gr$=683.985\mathrm{e}{3}$ (worst)}
       \end{overpic} 
 \caption{Relative errors over time for the velocity at Gr$=674.92\mathrm{e}{3}$ (left) and Gr$=683.985\mathrm{e}{3}$ (right). 
 }
 \label{Hess:batch2}
\end{center}
\end{figure}

The relative $L^2(\Omega)$ and $L^\infty(\Omega)$ errors over time for the best approximated case (Gr$=674.92\mathrm{e}{3}$) and the worst approximated case (Gr$=683.985\mathrm{e}{3}$) are shown in Fig.~\ref{Hess:batch2}. The main qualitative difference is that at Gr$=674.92\mathrm{e}{3}$ both the relative $L^2(\Omega)$ and $L^\infty(\Omega)$ errors oscillate around a fixed values, while at Gr$=683.985\mathrm{e}{3}$ they oscillate around a curved mean. 
An interesting feature of the errors at Gr$=683.985\mathrm{e}{3}$ is that the maximum error is after about $1000$ time steps and then the errors reduce again. See left panel in Fig.~\ref{Hess:batch2}. This is due to the fact that the phase of the approximation is most out-of-sync at about time step $1000$.

\begin{figure}[ht!]
\begin{center}
 \begin{overpic}[width=0.49\textwidth, grid = false]{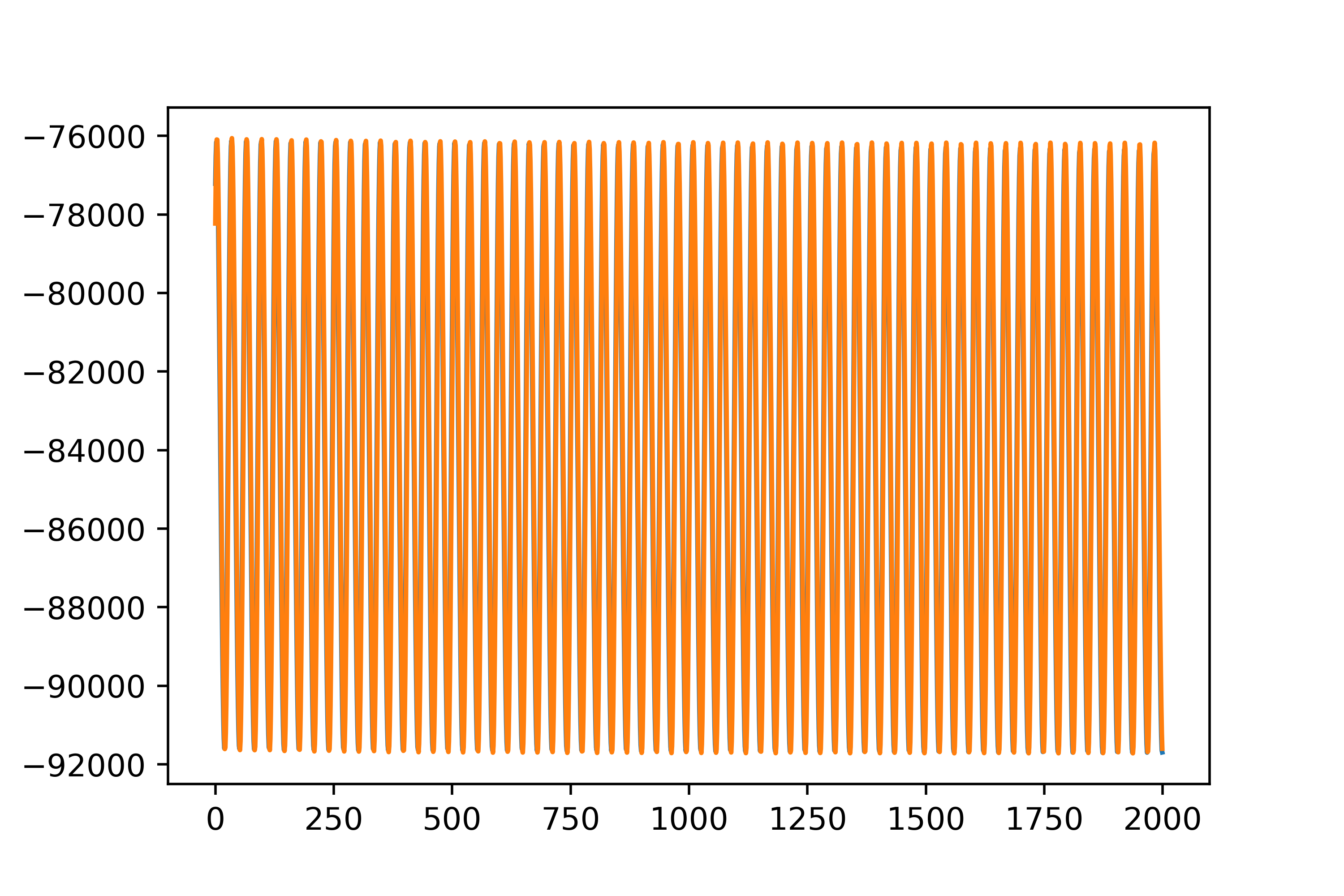}
        \put(25,60){horizontal component}
        \put(-10,35){\small{best}}
      \end{overpic} 
 \begin{overpic}[width=0.49\textwidth, grid = false]{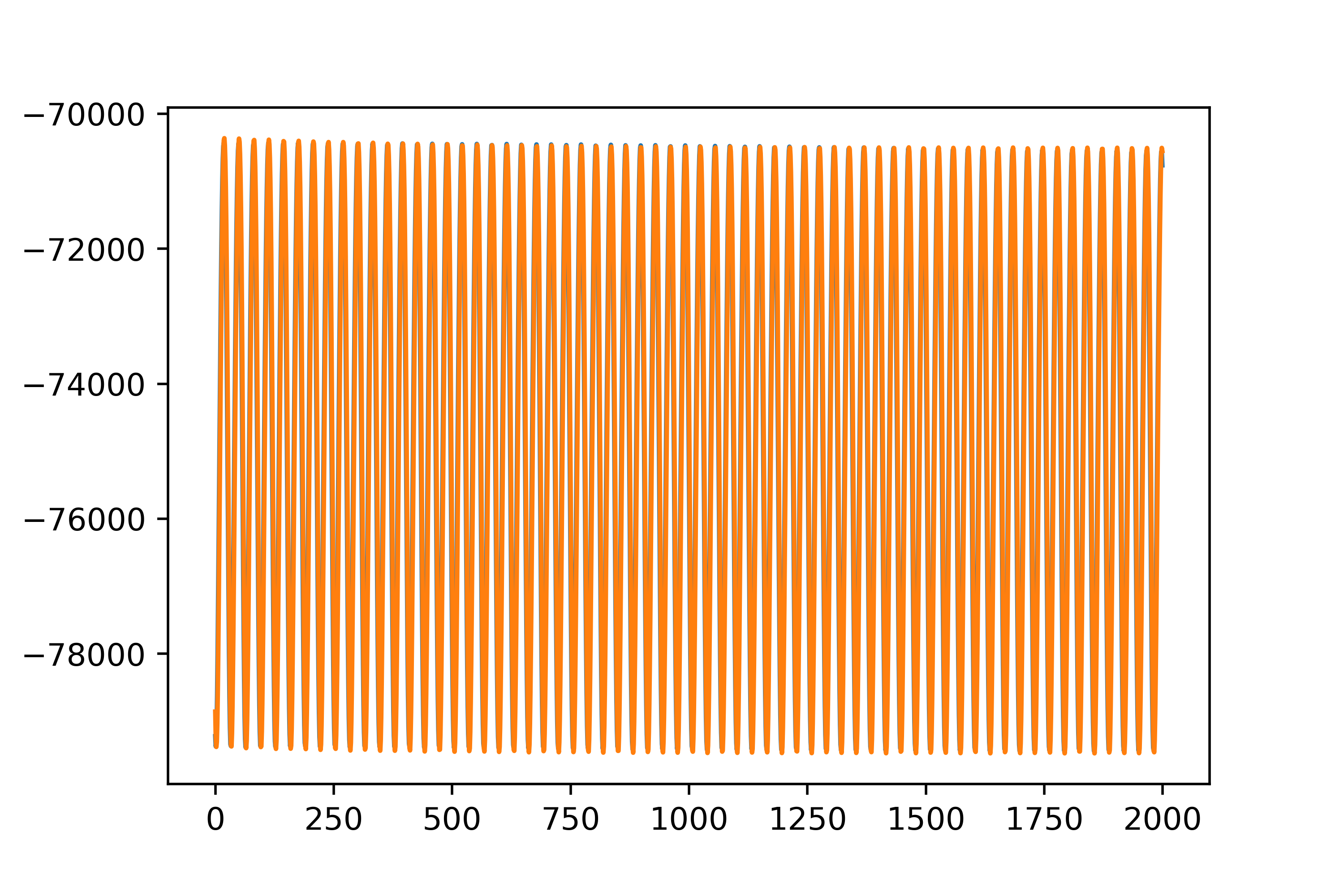} 
        \put(25,60){vertical component}
       \end{overpic} \\
   \begin{overpic}[width=0.49\textwidth, grid = false]{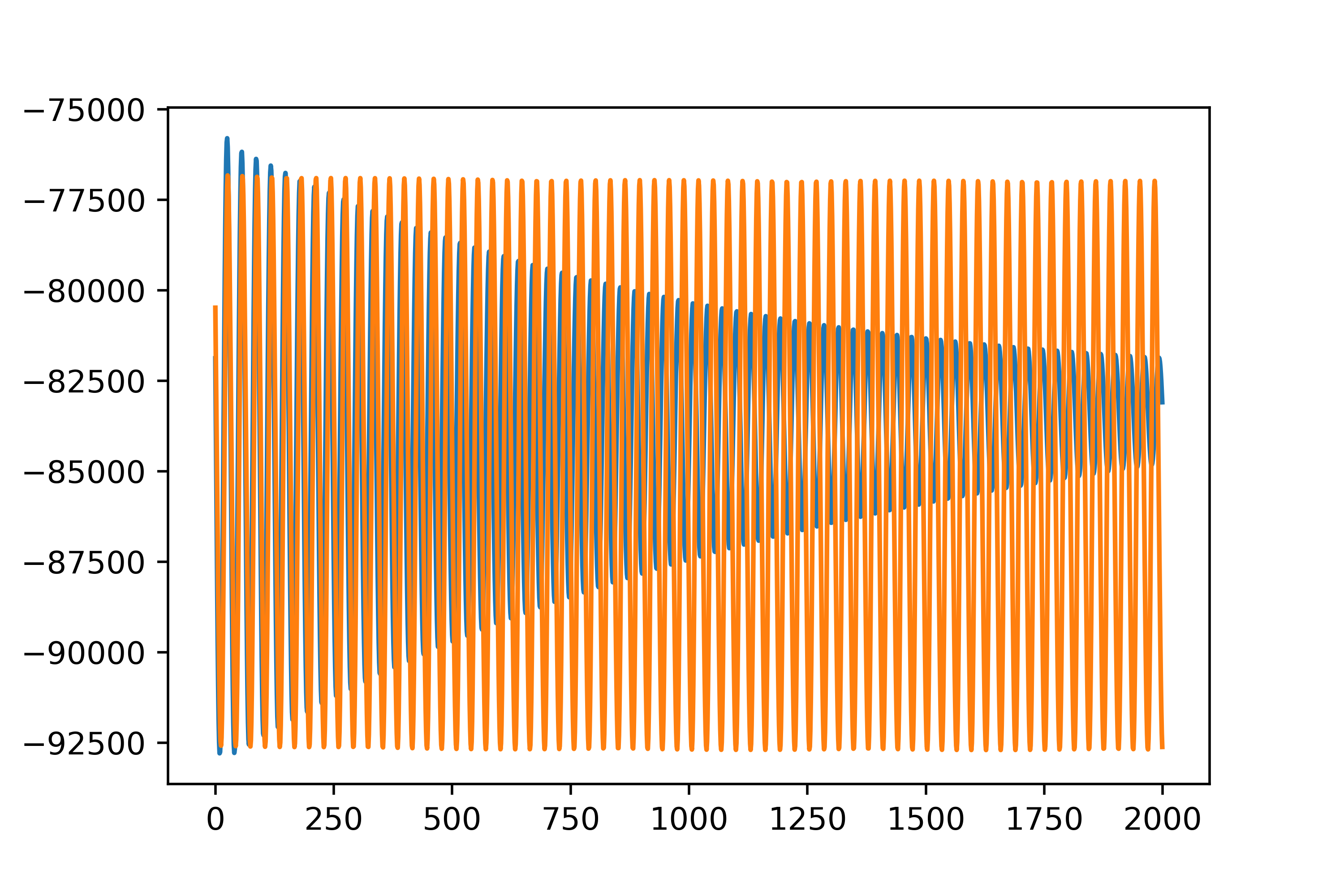}
       \put(-10,33){\small{worst}}
      \end{overpic} 
 \begin{overpic}[width=0.49\textwidth, grid = false]{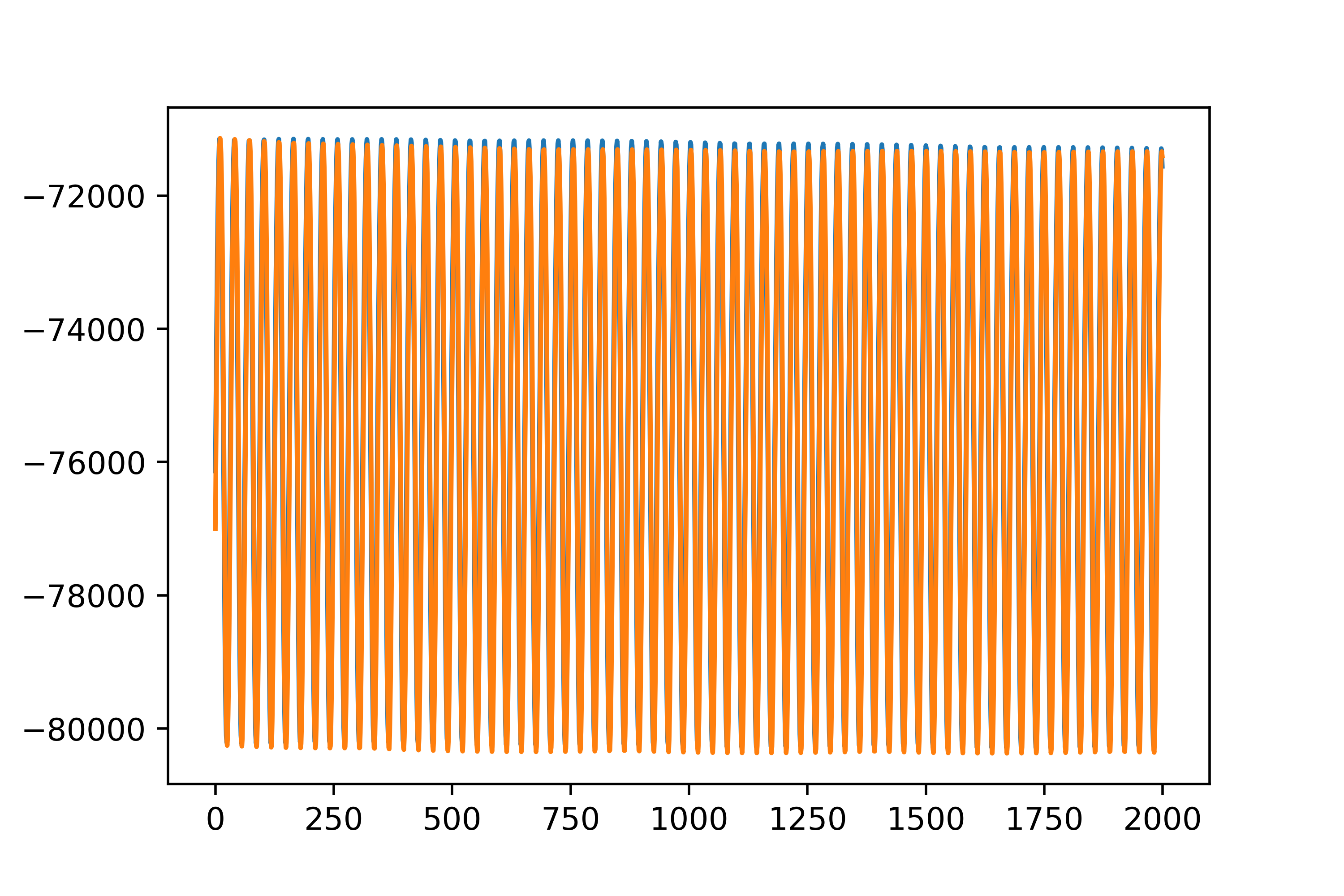} 
       \end{overpic} 
 \caption{Approximation of the first POD mode coefficient (blue) and reference solution (red) for the horizontal (left) and vertical (right) component of the velocity in the best approximation case, i.e.~Gr$=674.92\mathrm{e}{3}$ (top), and worst approximation case, i.e.~Gr$=683.985\mathrm{e}{3}$ (bottom).}
 \label{Hess:batch2_mode0_xy}
\end{center}
\end{figure}

\begin{figure}[ht!]
\begin{center}
 \includegraphics[width=0.49\textwidth]{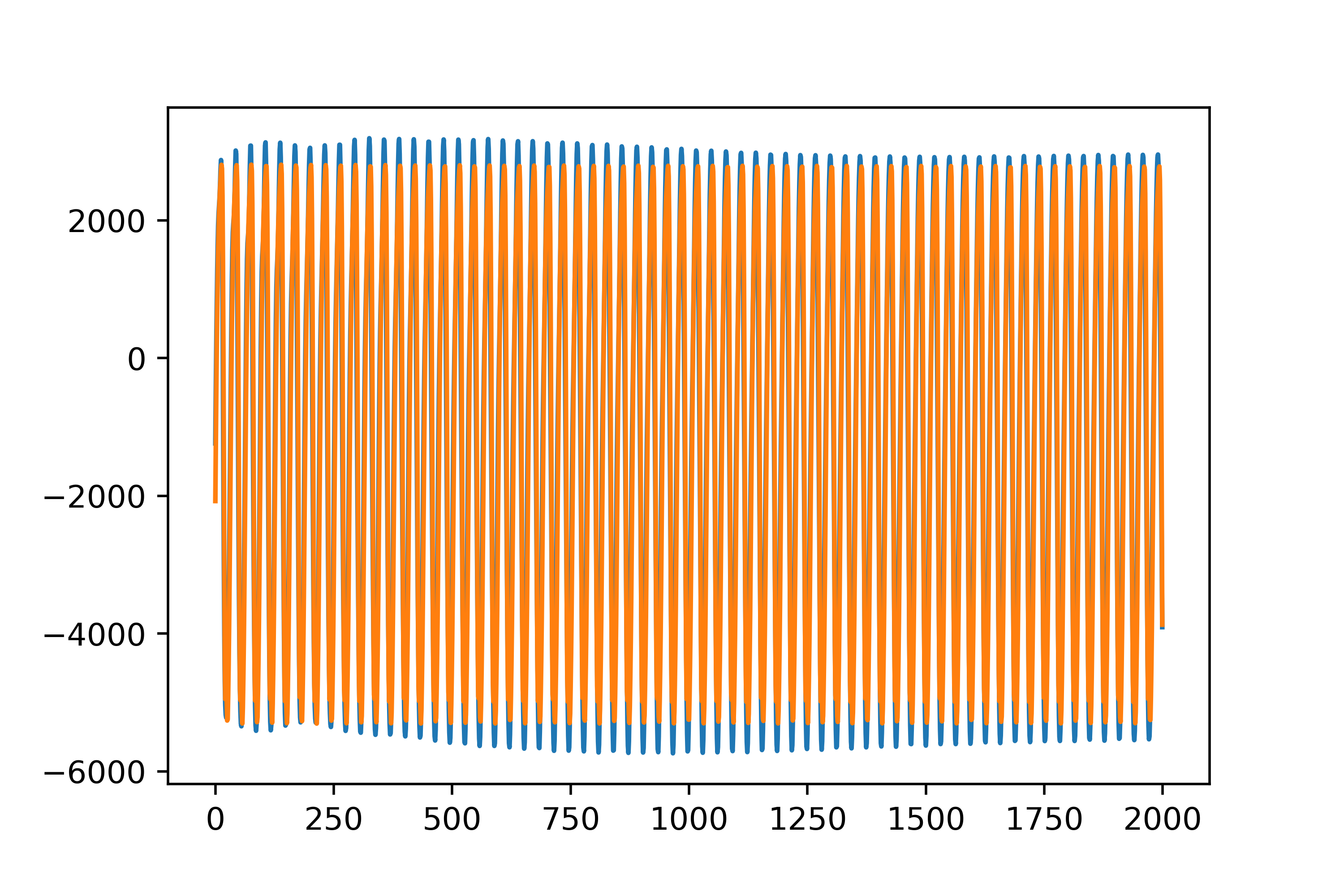} 
 \caption{Approximation of the fourth POD mode coefficient (blue) and reference solution (red) in the vertical velocity component for Gr$=674.92\mathrm{e}{3}$.}
 \label{Hess:batch2_rand6_mode3_y}
\end{center}
\end{figure}

Once again, it is instructive to look at how the first POD modes for the horizontal and vertical components of the velocity are resolved in both cases. See Fig.~\ref{Hess:batch2_mode0_xy}. 
We see that at Gr$=674.92\mathrm{e}{3}$ the blue and red curves are practically superimposed for both velocity components, indicating that error is negligible. Upon investigating the other POD modes, it becomes visible that the fourth mode for the vertical component of the velocity dominates the error. See Fig.~\ref{Hess:batch2_rand6_mode3_y}.
As for Gr$=683.985\mathrm{e}{3}$, the error originates from the approximation of the first mode in the horizontal component as shown in the bottom left panel of Fig.~\ref{Hess:batch2_mode0_xy}. 
Although this panel shows that the horizontal component of the velocity at Gr$=683.985\mathrm{e}{3}$ will not exhibit a 
time-periodic behaviour, we could not find a
number of POD modes and DMD modes to avoid this
while keeping the same training samples.
The accuracy could be improved by increasing the
number of training samples.

\begin{figure}[htb!]
\begin{center}
 \includegraphics[scale=.22]{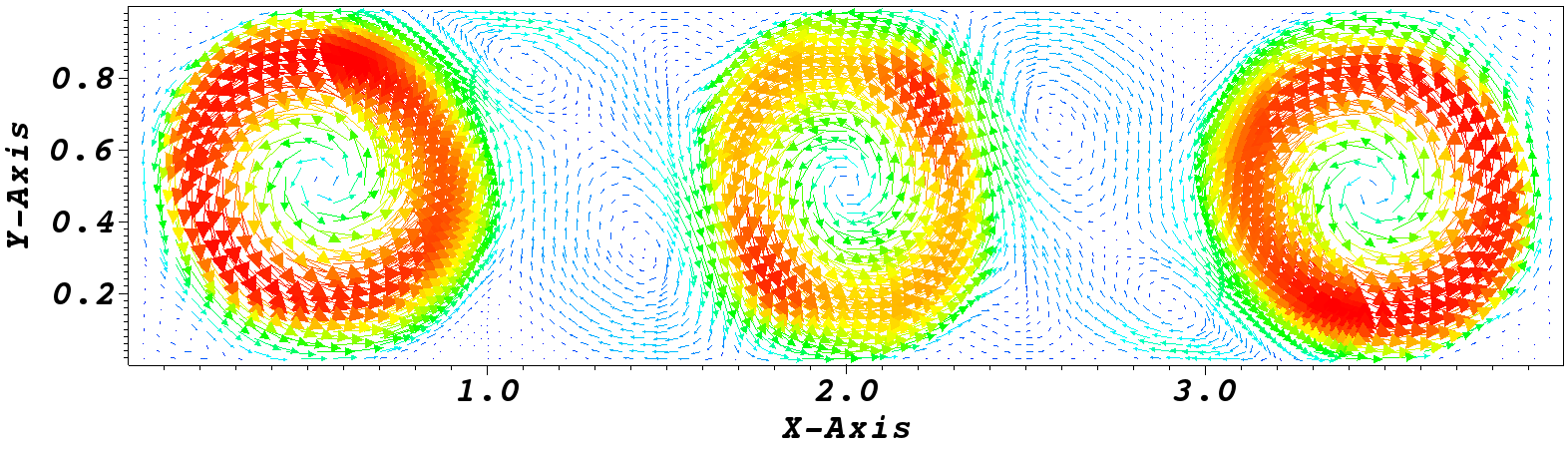} $\quad$
 \includegraphics[scale=.45]{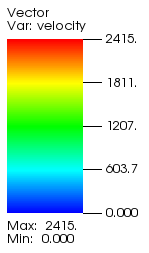} \\
 \includegraphics[scale=.22]{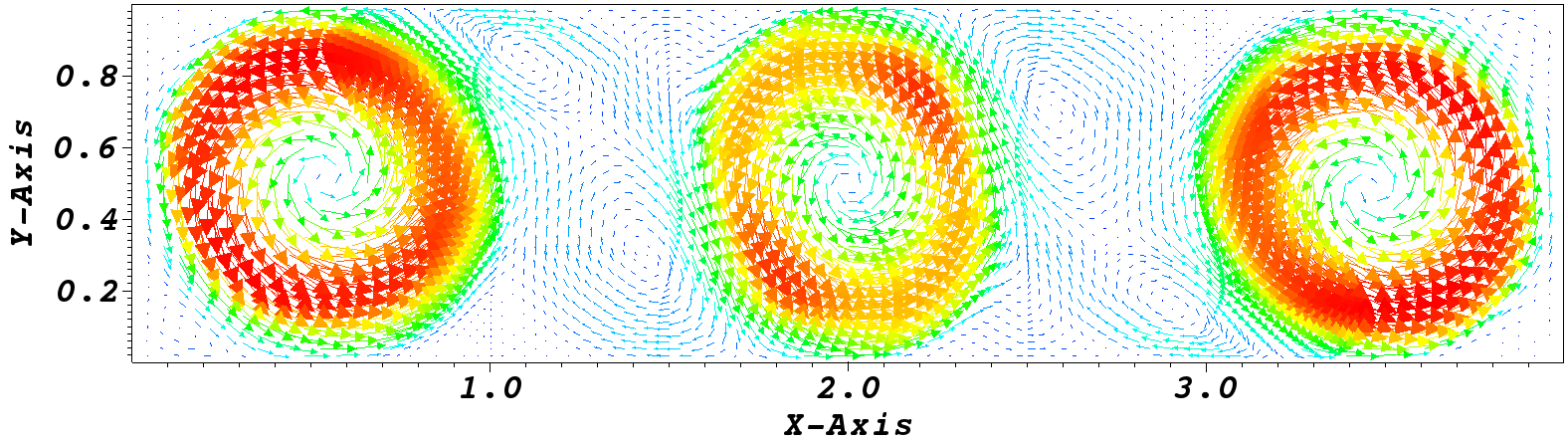} $\quad$
 \includegraphics[scale=.45]{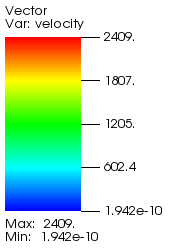} \\
 \includegraphics[scale=.22]{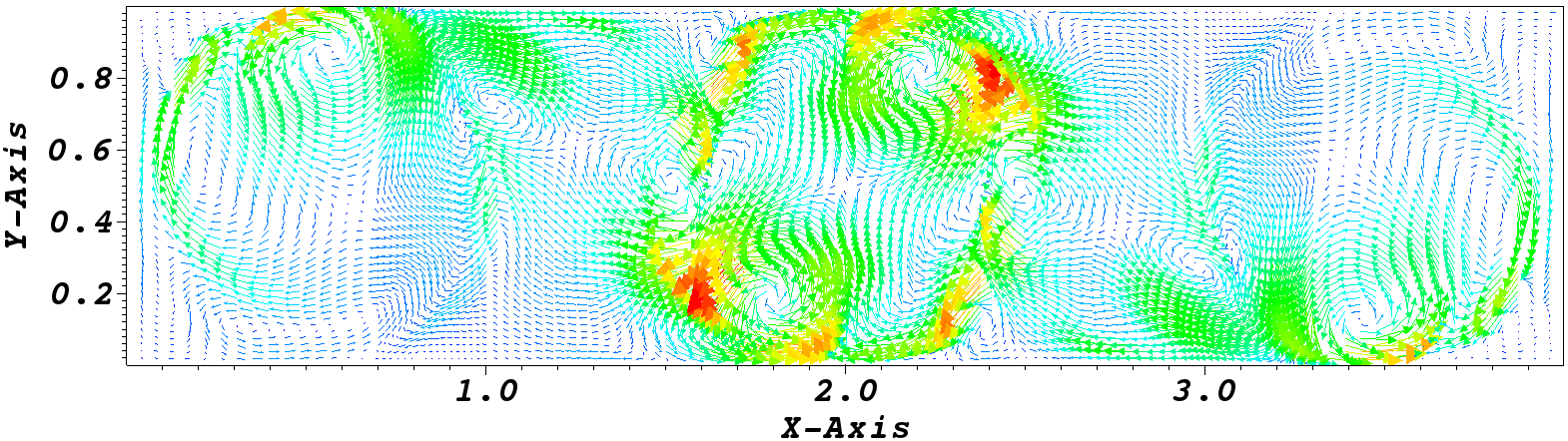} $\quad$
 \includegraphics[scale=.45]{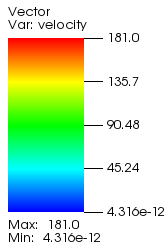} 
 \caption{Comparison of the velocity vector field obtained with the full order model (top) versus our ROM approach (middle) and the difference of the two (bottom) at Gr$=674.92\mathrm{e}{3}$ and after $13\mathrm{e}{3}$ time steps. 
 } 
 \label{Hess:FOM_vs_ROM_batch2_rand6}
\end{center}
\end{figure}

\begin{figure}[htb!]
\begin{center}
 \includegraphics[scale=.22]{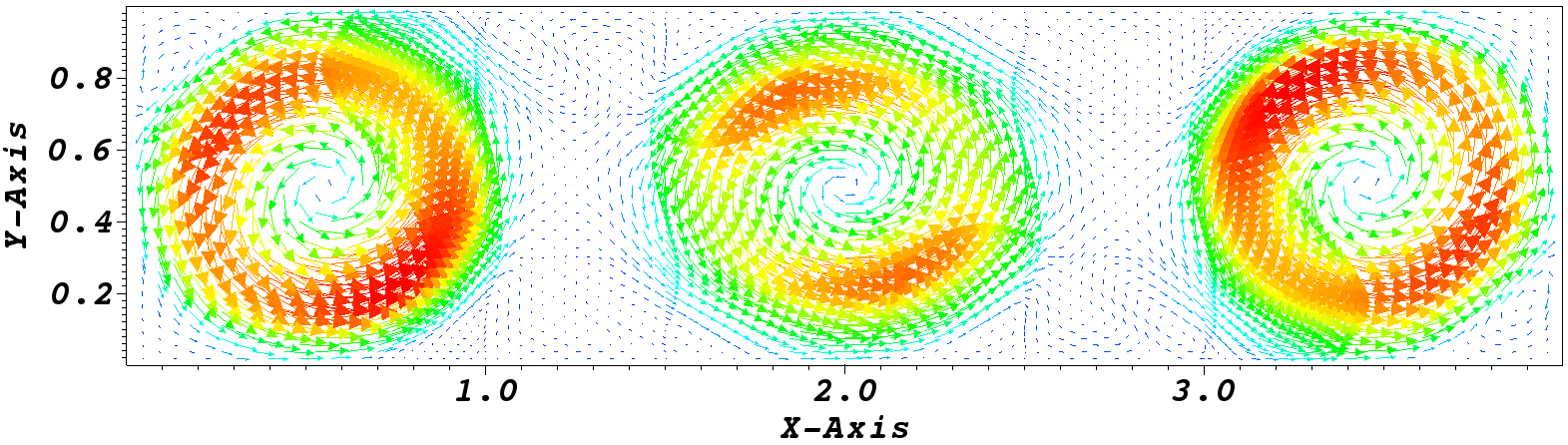} $\quad$
 \includegraphics[scale=.45]{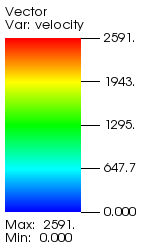} \\
 \includegraphics[scale=.22]{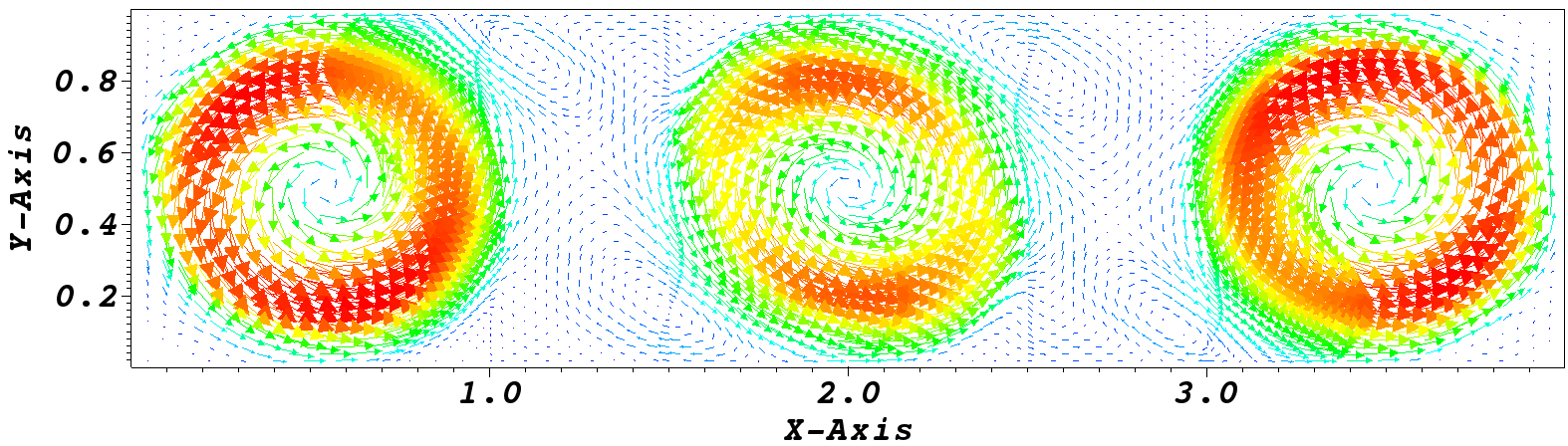} $\quad$
 \includegraphics[scale=.45]{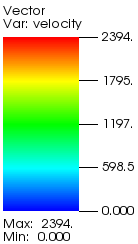} \\
 \includegraphics[scale=.22]{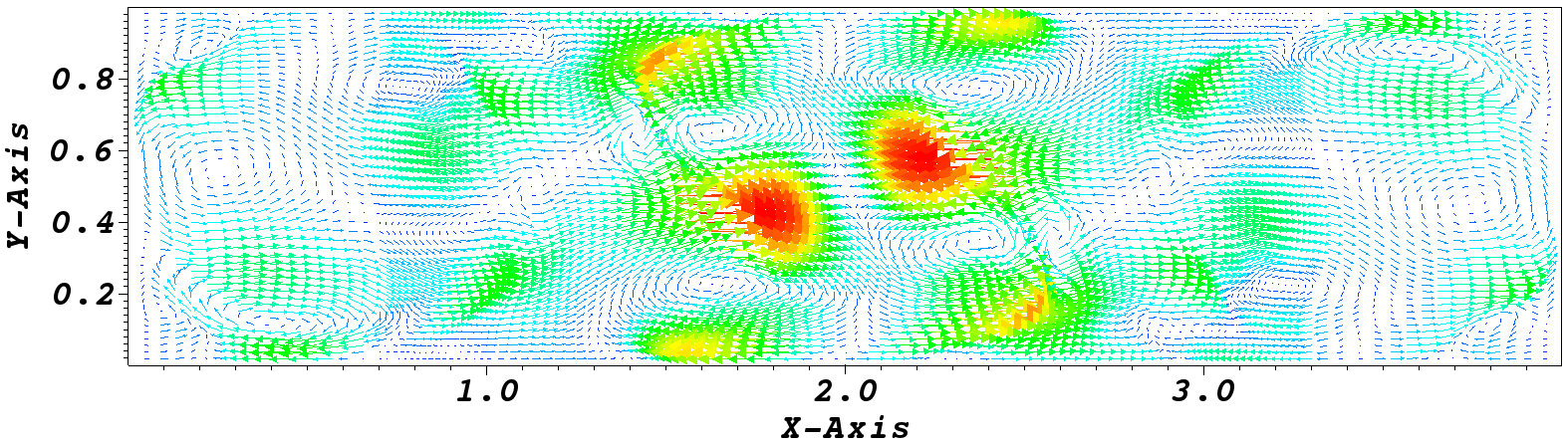} $\quad$
 \includegraphics[scale=.45]{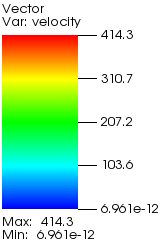} 
 \caption{Comparison of the velocity vector field obtained with the full order model (top) versus our ROM approach (middle) and the difference of the two (bottom) at Gr$=683.985\mathrm{e}{3}$ and after $13\mathrm{e}{3}$ time steps. 
}
 \label{Hess:FOM_vs_ROM_batch2_rand2}
\end{center}
\end{figure}

Next, we report a qualitative 
comparison between solutions obtained with the 
full order model and our ROM
approach at Gr$=674.92\mathrm{e}{3}$ (in Fig.~\ref{Hess:FOM_vs_ROM_batch2_rand6})
and Gr$=683.985\mathrm{e}{3}$ (in Fig.~\ref{Hess:FOM_vs_ROM_batch2_rand2}).
Although at Gr$=674.92\mathrm{e}{3}$ the mean and maximum errors in $L^2$ and $L^\infty$ norms  are between $3.5\%$ and $5\%$, Fig.~\ref{Hess:FOM_vs_ROM_batch2_rand6} shows that all important features have been captured by the ROM. On the other hand, Fig.~\ref{Hess:FOM_vs_ROM_batch2_rand2} shows that the middle roll given by our ROM is out of phase with respect to the one in the full order solution.
This could justify mean and maximum errors  in $L^2$ and $L^\infty$ norms of about $14\%$ and $25\%$, respectively, as reported in Table~\ref{high_Gr_results}.

Although the higher accuracy in the medium Grashof range can be attributed to less complex high-order simulations, some additional comments are in order.
Recall that in the $100\mathrm{e}{3} - 120\mathrm{e}{3}$ Grashof range the solution is not too far from a steady state. Comparing the solution at fixed time $t = 1e5$ produces $2\% - 13\%$ error, while our ROM method reduces this error to $2\%-4\%$, with the $2\%$ lower bound arising from the POD projection error.
In the high Grashof range the ROM error is $2\%-15\%$ on average, while comparing to the solution at fixed time $t = 1e5$  produces an error of $2\%-35\%$. Thus, relative to the mean field solution the model reduction works equally well in the medium and high Gr regimes.

\begin{figure}[ht!]
\begin{center}
 \begin{overpic}[width=0.49\textwidth, grid = false]{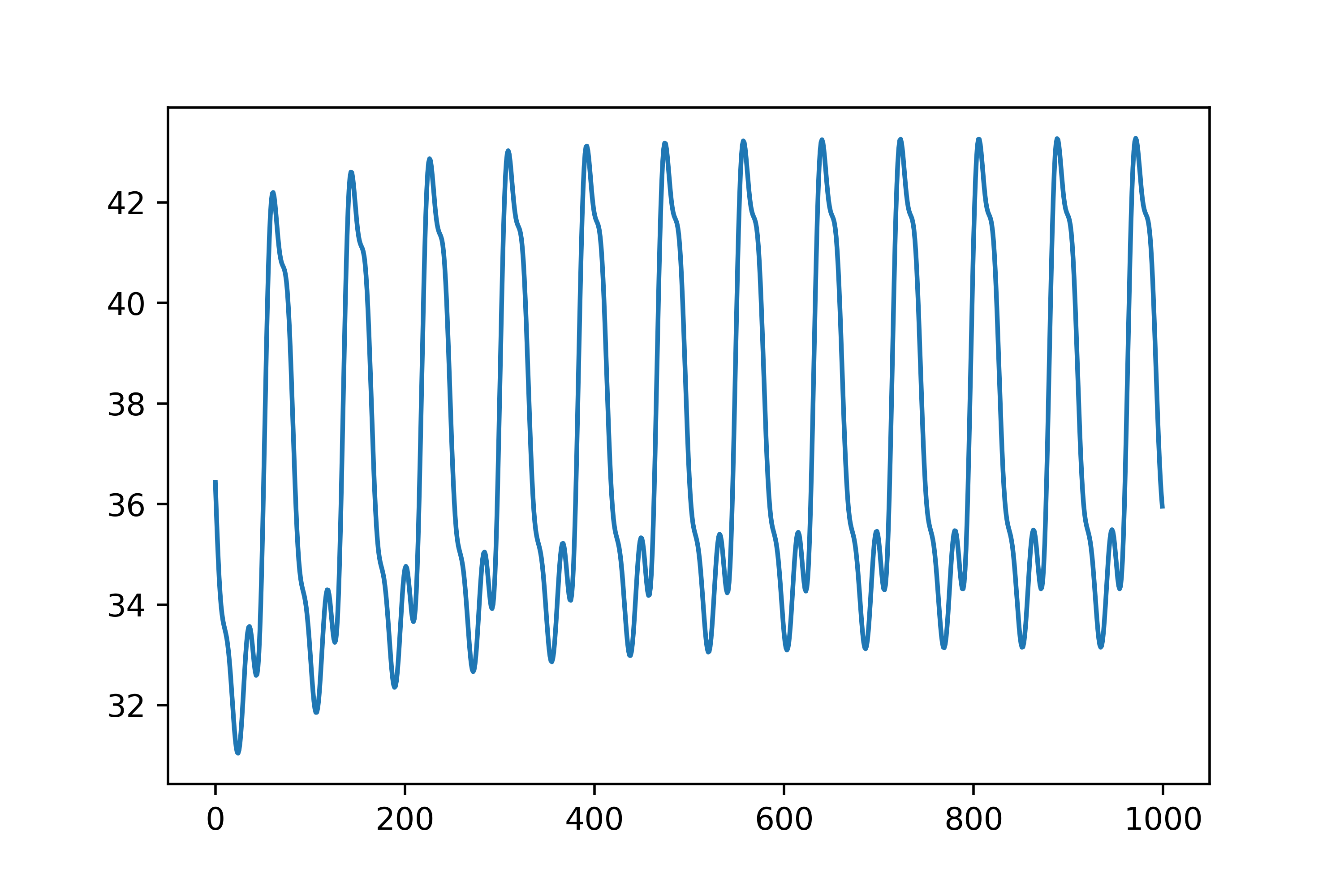} 
        \put(15,60){9-th mode at Gr$=113.8\mathrm{e}{3}$}
      \end{overpic} 
 \begin{overpic}[width=0.49\textwidth, grid = false]{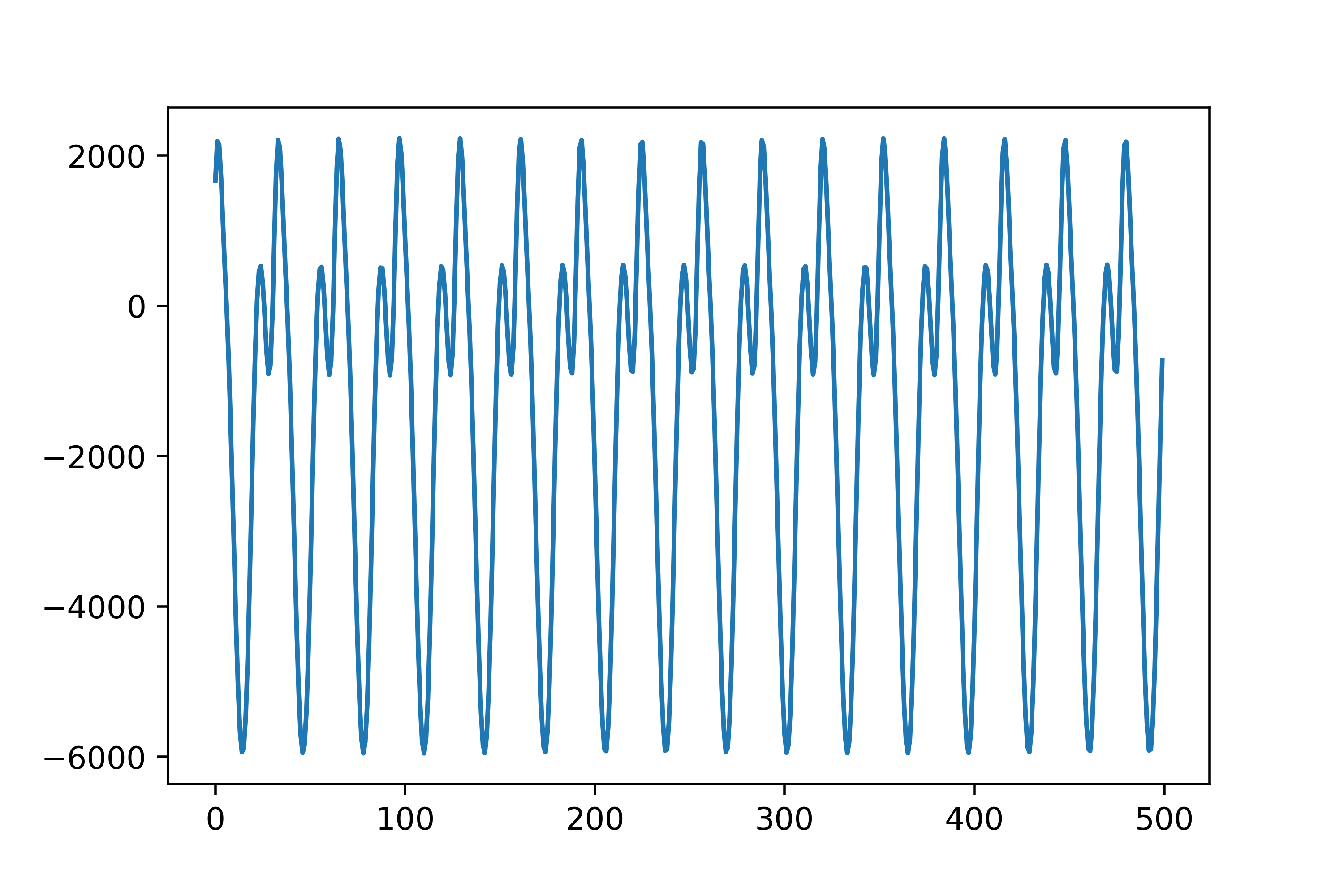} 
        \put(15,60){7-th mode at Gr$=663.8\mathrm{e}{3}$ }
       \end{overpic} 
 \caption{Evolution of the
 coefficient of the ninth mode at Gr$=113.8\mathrm{e}{3}$ in y-direction (left) and seventh mode at Gr$=663.8\mathrm{e}{3}$ in y-direction (right).}
 \label{Hess:higher_modes}
\end{center}
\end{figure}

We conclude this section by taking a look at the coefficients of the higher order modes. Fig.~\ref{Hess:higher_modes} shows that such coefficients can have more complex features than the 
the coefficients of the most dominant modes.
Since the coefficients of the first modes usually
have the largest amplitudes, these higher order modes are not well approximated by the DMD in general. However, it is more important to strive for a low approximation error in the first, most dominant modes than to accurately reproduce the higher order modes.

\section{Conclusions and future perspectives}

This work introduces a data-driven ROM approach to 
compute efficiently complex time-periodic simulations. There are three main building blocks: proper orthogonal decomposition, dynamic mode decomposition (DMD), and manifold interpolation. Our ROM approach is 
tested and validated on the Rayleigh-B\'{e}nard cavity problem with fixed aspect ratio and variable Grashof number (Gr). We focus on two parameter domains with time-periodic solutions: a medium Gr range and a high Gr range, which is close to turbulent behaviour. The key feature of our ROM is that it allows to recover frequencies not present in the sampled high-order solutions. This is crucial to achieve accurate simulations at new parameter values. Although in some instances of the high Gr regime, the mean relative error remained above $10\%$, most simulations achieved engineering accuracy. 

Our multi-stage ROM method could be further improved as follows. Stability of the DMD algorithm could be enforced by various techniques employed in the DMD literature. The manifold interpolation might benefit from using non-flat metrics to interpolate the reduced DMD operator and the use of the complex DMD could be explored. 
Finally, it would be interesting to apply the proposed approach to other practical engineering problems and higher dimensional parameter domains as well as quasi-periodic systems.

\bmhead{Acknowledgments}

We acknowledge the support provided by the European Research Council Executive Agency by the Consolidator Grant project AROMA-CFD ``Advanced Reduced Order Methods with Applications in Computational Fluid Dynamics" - GA 681447, H2020-ERC CoG 2015 AROMA-CFD, PI G. Rozza, and INdAM-GNCS 2019-2020 projects.
This work was also partially supported by US National Science Foundation through grant DMS-1953535. A.~Quaini acknowledges 
support from the Radcliffe Institute for Advanced Study at Harvard University where she has been a 2021-2022 William and Flora Hewlett Foundation Fellow.

\section*{Statements and Declarations}

\begin{itemize}
\item Conflict of interest/Competing interests 

The authors have no conflicts of interest or competing interests.

\end{itemize}

\bibliography{rbsissa,latexbi}


\begin{thebibliography}{51}
\ifx \bisbn   \undefined \def \bisbn  #1{ISBN #1}\fi
\ifx \binits  \undefined \def \binits#1{#1}\fi
\ifx \bauthor  \undefined \def \bauthor#1{#1}\fi
\ifx \batitle  \undefined \def \batitle#1{#1}\fi
\ifx \bjtitle  \undefined \def \bjtitle#1{#1}\fi
\ifx \bvolume  \undefined \def \bvolume#1{\textbf{#1}}\fi
\ifx \byear  \undefined \def \byear#1{#1}\fi
\ifx \bissue  \undefined \def \bissue#1{#1}\fi
\ifx \bfpage  \undefined \def \bfpage#1{#1}\fi
\ifx \blpage  \undefined \def \blpage #1{#1}\fi
\ifx \burl  \undefined \def \burl#1{\textsf{#1}}\fi
\ifx \doiurl  \undefined \def \doiurl#1{\url{https://doi.org/#1}}\fi
\ifx \betal  \undefined \def \betal{\textit{et al.}}\fi
\ifx \binstitute  \undefined \def \binstitute#1{#1}\fi
\ifx \binstitutionaled  \undefined \def \binstitutionaled#1{#1}\fi
\ifx \bctitle  \undefined \def \bctitle#1{#1}\fi
\ifx \beditor  \undefined \def \beditor#1{#1}\fi
\ifx \bpublisher  \undefined \def \bpublisher#1{#1}\fi
\ifx \bbtitle  \undefined \def \bbtitle#1{#1}\fi
\ifx \bedition  \undefined \def \bedition#1{#1}\fi
\ifx \bseriesno  \undefined \def \bseriesno#1{#1}\fi
\ifx \blocation  \undefined \def \blocation#1{#1}\fi
\ifx \bsertitle  \undefined \def \bsertitle#1{#1}\fi
\ifx \bsnm \undefined \def \bsnm#1{#1}\fi
\ifx \bsuffix \undefined \def \bsuffix#1{#1}\fi
\ifx \bparticle \undefined \def \bparticle#1{#1}\fi
\ifx \barticle \undefined \def \barticle#1{#1}\fi
\bibcommenthead
\ifx \bconfdate \undefined \def \bconfdate #1{#1}\fi
\ifx \botherref \undefined \def \botherref #1{#1}\fi
\ifx \url \undefined \def \url#1{\textsf{#1}}\fi
\ifx \bchapter \undefined \def \bchapter#1{#1}\fi
\ifx \bbook \undefined \def \bbook#1{#1}\fi
\ifx \bcomment \undefined \def \bcomment#1{#1}\fi
\ifx \oauthor \undefined \def \oauthor#1{#1}\fi
\ifx \citeauthoryear \undefined \def \citeauthoryear#1{#1}\fi
\ifx \endbibitem  \undefined \def \endbibitem {}\fi
\ifx \bconflocation  \undefined \def \bconflocation#1{#1}\fi
\ifx \arxivurl  \undefined \def \arxivurl#1{\textsf{#1}}\fi
\csname PreBibitemsHook\endcsname

\bibitem{HBMOR_vol1}
\begin{bbook}
\beditor{\bsnm{Benner}, \binits{P.}},
\beditor{\bsnm{Grivet-Talocia}, \binits{S.}},
\beditor{\bsnm{Quarteroni}, \binits{A.}},
\beditor{\bsnm{Rozza}, \binits{G.}},
\beditor{\bsnm{Schilders}, \binits{W.}},
\beditor{\bsnm{Silveira}, \binits{L.M.}} (eds.):
\bbtitle{Model Order Reduction: Volume 1: System- and Data-Driven Methods and
  Algorithms}.
\bpublisher{De Gruyter},
\blocation{Berlin, Boston}
(\byear{2021}).
\doiurl{10.1515/9783110498967}
\end{bbook}
\endbibitem

\bibitem{HBMOR_vol2}
\begin{bbook}
\beditor{\bsnm{Benner}, \binits{P.}},
\beditor{\bsnm{Grivet-Talocia}, \binits{S.}},
\beditor{\bsnm{Quarteroni}, \binits{A.}},
\beditor{\bsnm{Rozza}, \binits{G.}},
\beditor{\bsnm{Schilders}, \binits{W.}},
\beditor{\bsnm{Silveira}, \binits{L.M.}} (eds.):
\bbtitle{Model Order Reduction: Volume 2: Snapshot-Based Methods and
  Algorithms}.
\bpublisher{De Gruyter},
\blocation{Berlin, Boston}
(\byear{2021}).
\doiurl{10.1515/9783110671490}
\end{bbook}
\endbibitem

\bibitem{HBMOR_vol3}
\begin{bbook}
\beditor{\bsnm{Benner}, \binits{P.}},
\beditor{\bsnm{Grivet-Talocia}, \binits{S.}},
\beditor{\bsnm{Quarteroni}, \binits{A.}},
\beditor{\bsnm{Rozza}, \binits{G.}},
\beditor{\bsnm{Schilders}, \binits{W.}},
\beditor{\bsnm{Silveira}, \binits{L.M.}} (eds.):
\bbtitle{Model Order Reduction: Volume 3: Applications}.
\bpublisher{De Gruyter},
\blocation{Berlin, Boston}
(\byear{2021}).
\doiurl{10.1515/9783110499001}
\end{bbook}
\endbibitem

\bibitem{NOOR1982955}
\begin{barticle}
\bauthor{\bsnm{Noor}, \binits{A.}}:
\batitle{On making large nonlinear problems small}.
\bjtitle{Computer Methods in Applied Mechanics and Engineering}
\bvolume{34}(\bissue{1}),
\bfpage{955}--\blpage{985}
(\byear{1982})
\end{barticle}
\endbibitem

\bibitem{Noor:1994}
\begin{barticle}
\bauthor{\bsnm{Noor}, \binits{A.}}:
\batitle{Recent advances and applications of reduction methods}.
\bjtitle{ASME. Appl. Mech. Rev.}
\bvolume{5}(\bissue{47}),
\bfpage{125}--\blpage{146}
(\byear{1994})
\end{barticle}
\endbibitem

\bibitem{Noor:1983}
\begin{barticle}
\bauthor{\bsnm{Noor}, \binits{A.}},
\bauthor{\bsnm{Peters}, \binits{J.}}:
\batitle{Multiple-parameter reduced basis technique for bifurcation and
  post-buckling analyses of composite materiale}.
\bjtitle{International {J}ournal for {N}umerical {M}ethods in {E}ngineering}
\bvolume{19},
\bfpage{1783}--\blpage{1803}
(\byear{1983})
\end{barticle}
\endbibitem

\bibitem{NOOR198367}
\begin{barticle}
\bauthor{\bsnm{Noor}, \binits{A.}},
\bauthor{\bsnm{Peters}, \binits{J.}}:
\batitle{Recent advances in reduction methods for instability analysis of
  structures}.
\bjtitle{Computers \& Structures}
\bvolume{16}(\bissue{1}),
\bfpage{67}--\blpage{80}
(\byear{1983})
\end{barticle}
\endbibitem

\bibitem{Maday:RB2}
\begin{barticle}
\bauthor{\bsnm{Herrero}, \binits{H.}},
\bauthor{\bsnm{Maday}, \binits{Y.}},
\bauthor{\bsnm{Pla}, \binits{F.}}:
\batitle{{RB (Reduced Basis) for RB (Rayleigh-B\'enard)}}.
\bjtitle{Computer Methods in Applied Mechanics and Engineering}
\bvolume{261-262},
\bfpage{132}--\blpage{141}
(\byear{2013})
\end{barticle}
\endbibitem

\bibitem{PLA2015162}
\begin{barticle}
\bauthor{\bsnm{Pla}, \binits{F.}},
\bauthor{\bsnm{Herrero}, \binits{H.}},
\bauthor{\bsnm{Vega}, \binits{J.}}:
\batitle{A flexible symmetry-preserving {G}alerkin/{POD} reduced order model
  applied to a convective instability problem}.
\bjtitle{Computers \& Fluids}
\bvolume{119},
\bfpage{162}--\blpage{175}
(\byear{2015})
\end{barticle}
\endbibitem

\bibitem{PR15}
\begin{barticle}
\bauthor{\bsnm{Pitton}, \binits{G.}},
\bauthor{\bsnm{Rozza}, \binits{G.}}:
\batitle{On the application of reduced basis methods to bifurcation problems in
  incompressible fluid dynamics}.
\bjtitle{Journal of Scientific Computing}
\bvolume{73}(\bissue{1}),
\bfpage{157}--\blpage{177}
(\byear{2017})
\end{barticle}
\endbibitem

\bibitem{PITTON2017534}
\begin{barticle}
\bauthor{\bsnm{Pitton}, \binits{G.}},
\bauthor{\bsnm{Quaini}, \binits{A.}},
\bauthor{\bsnm{Rozza}, \binits{G.}}:
\batitle{Computational reduction strategies for the detection of steady
  bifurcations in incompressible fluid-dynamics: Applications to {C}oanda
  effect in cardiology}.
\bjtitle{Journal of Computational Physics}
\bvolume{344},
\bfpage{534}--\blpage{557}
(\byear{2017})
\end{barticle}
\endbibitem

\bibitem{PichiRozza2019}
\begin{barticle}
\bauthor{\bsnm{Pichi}, \binits{F.}},
\bauthor{\bsnm{Rozza}, \binits{G.}}:
\batitle{Reduced basis approaches for parametrized bifurcation problems held by
  non-linear von {K}{\'a}rm{\'a}n equations}.
\bjtitle{Journal of Scientific Computing}
\bvolume{81},
\bfpage{112}--\blpage{135}
(\byear{2019}).
\doiurl{10.1007/s10915-019-01003-3}
\end{barticle}
\endbibitem

\bibitem{pichi2020optcntrl}
\begin{botherref}
\oauthor{\bsnm{Pichi}, \binits{F.}},
\oauthor{\bsnm{Strazzullo}, \binits{M.}},
\oauthor{\bsnm{Ballarin}, \binits{F.}},
\oauthor{\bsnm{Rozza}, \binits{G.}}:
Driving bifurcating parametrized nonlinear pdes by optimal control strategies:
  application to navier-stokes equations with model order reduction.
ArXiv preprint
(2020)
\end{botherref}
\endbibitem

\bibitem{pichi2021fsi}
\begin{botherref}
\oauthor{\bsnm{Khamlich}, \binits{M.}},
\oauthor{\bsnm{Pichi}, \binits{F.}},
\oauthor{\bsnm{Rozza}, \binits{G.}}:
Model order reduction for bifurcating phenomena in fluid-structure interaction
  problems.
ArXiv preprint
(2021)
\end{botherref}
\endbibitem

\bibitem{PichiQuainiRozza}
\begin{barticle}
\bauthor{\bsnm{Pichi}, \binits{F.}},
\bauthor{\bsnm{Quaini}, \binits{A.}},
\bauthor{\bsnm{Rozza}, \binits{G.}}:
\batitle{A reduced order modeling technique to study bifurcating phenomena:
  Application to the {G}ross--{P}itaevskii equation}.
\bjtitle{SIAM Journal on Scientific Computing}
\bvolume{42}(\bissue{5}),
\bfpage{1115}--\blpage{1135}
(\byear{2020})
\end{barticle}
\endbibitem

\bibitem{BTBK14}
\begin{barticle}
\bauthor{\bsnm{Brunton}, \binits{S.L.}},
\bauthor{\bsnm{Tu}, \binits{J.H.}},
\bauthor{\bsnm{Bright}, \binits{I.}},
\bauthor{\bsnm{Kutz}, \binits{J.N.}}:
\batitle{Compressive sensing and low-rank libraries for classification of
  bifurcation regimes in nonlinear dynamical systems}.
\bjtitle{SIAM Journal on Applied Dynamical Systems}
\bvolume{13}(\bissue{4}),
\bfpage{1716}--\blpage{1732}
(\byear{2014}).
\doiurl{10.1137/130949282}
\end{barticle}
\endbibitem

\bibitem{KGBNB17}
\begin{barticle}
\bauthor{\bsnm{Kramer}, \binits{B.}},
\bauthor{\bsnm{Grover}, \binits{P.}},
\bauthor{\bsnm{Boufounos}, \binits{P.}},
\bauthor{\bsnm{Nabi}, \binits{S.}},
\bauthor{\bsnm{Benosman}, \binits{M.}}:
\batitle{Sparse sensing and {DMD}-based identification of flow regimes and
  bifurcations in complex flows}.
\bjtitle{SIAM Journal on Applied Dynamical Systems}
\bvolume{16}(\bissue{2}),
\bfpage{1164}--\blpage{1196}
(\byear{2017}).
\doiurl{10.1137/15M104565X}
\end{barticle}
\endbibitem

\bibitem{HessQuainiRozza2022_ETNA}
\begin{barticle}
\bauthor{\bsnm{Hess}, \binits{M.W.}},
\bauthor{\bsnm{Quaini}, \binits{A.}},
\bauthor{\bsnm{Rozza}, \binits{G.}}:
\batitle{A comparison of reduced-order modeling approaches using artificial
  neural networks for {PDE}s with bifurcating solutions}.
\bjtitle{ETNA - Electronic Transactions on Numerical Analysis}
\bvolume{56},
\bfpage{52}--\blpage{65}
(\byear{2022}).
\doiurl{10.1553/etna_vol56s52}
\end{barticle}
\endbibitem

\bibitem{Hess2019CMAME}
\begin{barticle}
\bauthor{\bsnm{Hess}, \binits{M.W.}},
\bauthor{\bsnm{Alla}, \binits{A.}},
\bauthor{\bsnm{Quaini}, \binits{A.}},
\bauthor{\bsnm{Rozza}, \binits{G.}},
\bauthor{\bsnm{Gunzburger}, \binits{M.}}:
\batitle{A localized reduced-order modeling approach for {PDE}s with
  bifurcating solutions}.
\bjtitle{Comput. Methods Appl. Mech. Engrg.}
\bvolume{351},
\bfpage{379}--\blpage{403}
(\byear{2019})
\end{barticle}
\endbibitem

\bibitem{POD_NN_Hesthaven_Ubbiali}
\begin{botherref}
\oauthor{\bsnm{Hesthaven}, \binits{J.}},
\oauthor{\bsnm{Ubbiali}, \binits{S.}}:
Non-intrusive reduced order modeling of nonlinear problems using neural
  networks.
Journal of Computational Physics
\textbf{363}
(2018).
\doiurl{10.1016/j.jcp.2018.02.037}
\end{botherref}
\endbibitem

\bibitem{pichi2021artificial}
\begin{botherref}
\oauthor{\bsnm{Pichi}, \binits{F.}},
\oauthor{\bsnm{Ballarin}, \binits{F.}},
\oauthor{\bsnm{Rozza}, \binits{G.}},
\oauthor{\bsnm{Hesthaven}, \binits{J.S.}}:
An artificial neural network approach to bifurcating phenomena in computational
  fluid dynamics.
ArXiv preprint
(2021)
\end{botherref}
\endbibitem

\bibitem{10.5555/3327757.3327764}
\begin{bchapter}
\bauthor{\bsnm{Chen}, \binits{R.T.Q.}},
\bauthor{\bsnm{Rubanova}, \binits{Y.}},
\bauthor{\bsnm{Bettencourt}, \binits{J.}},
\bauthor{\bsnm{Duvenaud}, \binits{D.}}:
\bctitle{Neural ordinary differential equations}.
In: \bbtitle{Proceedings of the 32nd International Conference on Neural
  Information Processing Systems}.
\bsertitle{NIPS'18},
pp. \bfpage{6572}--\blpage{6583}.
\bpublisher{Curran Associates Inc.},
\blocation{Red Hook, NY, USA}
(\byear{2018})
\end{bchapter}
\endbibitem

\bibitem{Brunton3932}
\begin{barticle}
\bauthor{\bsnm{Brunton}, \binits{S.L.}},
\bauthor{\bsnm{Proctor}, \binits{J.L.}},
\bauthor{\bsnm{Kutz}, \binits{J.N.}}:
\batitle{Discovering governing equations from data by sparse identification of
  nonlinear dynamical systems}.
\bjtitle{Proceedings of the National Academy of Sciences}
\bvolume{113}(\bissue{15}),
\bfpage{3932}--\blpage{3937}
(\byear{2016}).
\doiurl{10.1073/pnas.1517384113}
\end{barticle}
\endbibitem

\bibitem{Koopman1931}
\begin{barticle}
\bauthor{\bsnm{Koopman}, \binits{B.O.}}:
\batitle{Hamiltonian systems and transformation in {Hilbert} space}.
\bjtitle{Proceedings of the National Academy of Sciences of the United States
  of America}
\bvolume{17}(\bissue{5}),
\bfpage{315}--\blpage{318}
(\byear{1931}).
\doiurl{10.1073/pnas.17.5.315}
\end{barticle}
\endbibitem

\bibitem{doi:10.1137/1.9781611974508}
\begin{bbook}
\bauthor{\bsnm{Kutz}, \binits{J.N.}},
\bauthor{\bsnm{Brunton}, \binits{S.L.}},
\bauthor{\bsnm{Brunton}, \binits{B.W.}},
\bauthor{\bsnm{Proctor}, \binits{J.L.}}:
\bbtitle{Dynamic Mode Decomposition}.
\bpublisher{Society for Industrial and Applied Mathematics},
\blocation{Philadelphia, PA}
(\byear{2016}).
\doiurl{10.1137/1.9781611974508}
\end{bbook}
\endbibitem

\bibitem{schmid_2010}
\begin{barticle}
\bauthor{\bsnm{Schmid}, \binits{P.J.}}:
\batitle{Dynamic mode decomposition of numerical and experimental data}.
\bjtitle{Journal of Fluid Mechanics}
\bvolume{656},
\bfpage{5}--\blpage{28}
(\byear{2010}).
\doiurl{10.1017/S0022112010001217}
\end{barticle}
\endbibitem

\bibitem{GAO2021110907}
\begin{botherref}
\oauthor{\bsnm{Gao}, \binits{Z.}},
\oauthor{\bsnm{Lin}, \binits{Y.}},
\oauthor{\bsnm{Sun}, \binits{X.}},
\oauthor{\bsnm{Zeng}, \binits{X.}}:
A reduced order method for nonlinear parameterized partial differential
  equations using dynamic mode decomposition coupled with k-nearest-neighbors
  regression.
Journal of Computational Physics,
110907
(2021).
\doiurl{10.1016/j.jcp.2021.110907}
\end{botherref}
\endbibitem

\bibitem{doi:10.1063/1.4913868}
\begin{barticle}
\bauthor{\bsnm{Sayadi}, \binits{T.}},
\bauthor{\bsnm{Schmid}, \binits{P.J.}},
\bauthor{\bsnm{Richecoeur}, \binits{F.}},
\bauthor{\bsnm{Durox}, \binits{D.}}:
\batitle{Parametrized data-driven decomposition for bifurcation analysis, with
  application to thermo-acoustically unstable systems}.
\bjtitle{Physics of Fluids}
\bvolume{27}(\bissue{3}),
\bfpage{037102}
(\byear{2015}).
\doiurl{10.1063/1.4913868}
\end{barticle}
\endbibitem

\bibitem{TezzeleDemoStabileMolaRozza2020}
\begin{botherref}
\oauthor{\bsnm{Tezzele}, \binits{M.}},
\oauthor{\bsnm{Demo}, \binits{N.}},
\oauthor{\bsnm{Stabile}, \binits{G.}},
\oauthor{\bsnm{Mola}, \binits{A.}},
\oauthor{\bsnm{Rozza}, \binits{G.}}:
{Enhancing CFD predictions in shape design problems by model and parameter
  space reduction}.
Advanced Modeling and Simulation in Engineering Sciences
\textbf{7}(40)
(2020).
\doiurl{10.1186/s40323-020-00177-y}
\end{botherref}
\endbibitem

\bibitem{andreuzzi2021dynamic}
\begin{botherref}
\oauthor{\bsnm{Andreuzzi}, \binits{F.}},
\oauthor{\bsnm{Demo}, \binits{N.}},
\oauthor{\bsnm{Rozza}, \binits{G.}}:
A dynamic mode decomposition extension for the forecasting of parametric
  dynamical systems.
ArXiv preprint
(2021)
\end{botherref}
\endbibitem

\bibitem{Zimmermann2019ManifoldIA}
\begin{bchapter}
\bauthor{\bsnm{Zimmermann}, \binits{R.}}:
\bctitle{Manifold interpolation}.
In: \bbtitle{Volume 1 System- and Data-Driven Methods and Algorithms},
pp. \bfpage{229}--\blpage{274}.
\bpublisher{De Gruyter},
\blocation{Berlin, Boston}
(\byear{2021}).
\doiurl{10.1515/9783110498967-007}
\end{bchapter}
\endbibitem

\bibitem{Roux:GAMM}
\begin{bbook}
\beditor{\bsnm{Roux}, \binits{B.}} (ed.):
\bbtitle{Numerical {S}imulation of {O}scillatory {C}onvection in {Low-Pr
  F}luids}.
\bsertitle{Notes on {N}umerical {F}luid {M}echanics and {M}ultidisciplinary
  {D}esign},
vol. \bseriesno{27}.
\bpublisher{Springer},
\blocation{Vieweg+Teubner Verlag}
(\byear{1990})
\end{bbook}
\endbibitem

\bibitem{Gelfgat:Ref11}
\begin{barticle}
\bauthor{\bsnm{Gelfgat}, \binits{A.Y.}},
\bauthor{\bsnm{Bar-Yoseph}, \binits{P.Z.}},
\bauthor{\bsnm{Yarin}, \binits{A.L.}}:
\batitle{Stability of multiple steady states of convection in laterally heated
  cavities}.
\bjtitle{Journal of {F}luid {M}echanics}
\bvolume{388},
\bfpage{315}--\blpage{334}
(\byear{1999})
\end{barticle}
\endbibitem

\bibitem{KAKIMOTO1995191}
\begin{barticle}
\bauthor{\bsnm{Kakimoto}, \binits{K.}}:
\batitle{Flow instability during crystal growth from the melt}.
\bjtitle{Progress in Crystal Growth and Characterization of Materials}
\bvolume{30}(\bissue{2}),
\bfpage{191}--\blpage{215}
(\byear{1995}).
\doiurl{10.1016/0960-8974(94)00013-J}
\end{barticle}
\endbibitem

\bibitem{Guermond_Shen_VCS}
\begin{barticle}
\bauthor{\bsnm{Guermond}, \binits{J.L.}},
\bauthor{\bsnm{Shen}, \binits{J.}}:
\batitle{Velocity-correction projection methods for incompressible flows}.
\bjtitle{SIAM Journal on Numerical Analysis}
\bvolume{41}(\bissue{1}),
\bfpage{112}--\blpage{134}
(\byear{2003})
\end{barticle}
\endbibitem

\bibitem{Karniadakis_Orszag_Israeli_splitting_methods}
\begin{barticle}
\bauthor{\bsnm{{Karniadakis}}, \binits{G.E.}},
\bauthor{\bsnm{{Orszag}}, \binits{S.A.}},
\bauthor{\bsnm{{Israeli}}, \binits{M.}}:
\batitle{{High-order splitting methods for the incompressible Navier-Stokes
  equations}}.
\bjtitle{Journal of Computational Physics}
\bvolume{97},
\bfpage{414}--\blpage{443}
(\byear{1991})
\end{barticle}
\endbibitem

\bibitem{hesthaven2015certified}
\begin{bbook}
\bauthor{\bsnm{Hesthaven}, \binits{J.}},
\bauthor{\bsnm{Rozza}, \binits{G.}},
\bauthor{\bsnm{Stamm}, \binits{B.}}:
\bbtitle{Certified Reduced Basis Methods for Parametrized Partial Differential
  Equations}.
\bpublisher{Springer}, \blocation{???}
(\byear{2015})
\end{bbook}
\endbibitem

\bibitem{LMQR:2014}
\begin{bchapter}
\bauthor{\bsnm{Lassila}, \binits{T.}},
\bauthor{\bsnm{Manzoni}, \binits{A.}},
\bauthor{\bsnm{Quarteroni}, \binits{A.}},
\bauthor{\bsnm{Rozza}, \binits{G.}}:
\bctitle{Model order reduction in fluid dynamics: challenges and perspectives}.
In: \beditor{\bsnm{Quarteroni}, \binits{A.}},
\beditor{\bsnm{Rozza}, \binits{G.}} (eds.)
\bbtitle{Reduced {O}rder {M}ethods for Modeling and Computational Reduction}.
\bsertitle{Modeling, {S}imulation and {A}pplications},
vol. \bseriesno{9},
pp. \bfpage{235}--\blpage{273}.
\bpublisher{Springer},
\blocation{Milano}
(\byear{2014}).
\bcomment{Chap. 9}
\end{bchapter}
\endbibitem

\bibitem{10.1007/978-3-319-10705-9_42}
\begin{bchapter}
\bauthor{\bsnm{Benner}, \binits{P.}},
\bauthor{\bsnm{Feng}, \binits{L.}},
\bauthor{\bsnm{Li}, \binits{S.}},
\bauthor{\bsnm{Zhang}, \binits{Y.}}:
\bctitle{Reduced-order modeling and rom-based optimization of batch
  chromatography}.
In: \beditor{\bsnm{Abdulle}, \binits{A.}},
\beditor{\bsnm{Deparis}, \binits{S.}},
\beditor{\bsnm{Kressner}, \binits{D.}},
\beditor{\bsnm{Nobile}, \binits{F.}},
\beditor{\bsnm{Picasso}, \binits{M.}} (eds.)
\bbtitle{Numerical Mathematics and Advanced Applications - ENUMATH 2013},
pp. \bfpage{427}--\blpage{435}.
\bpublisher{Springer},
\blocation{Cham}
(\byear{2015})
\end{bchapter}
\endbibitem

\bibitem{Demo2018}
\begin{barticle}
\bauthor{\bsnm{Demo}, \binits{N.}},
\bauthor{\bsnm{Tezzele}, \binits{M.}},
\bauthor{\bsnm{Rozza}, \binits{G.}}:
\batitle{Pydmd: Python dynamic mode decomposition}.
\bjtitle{Journal of Open Source Software}
\bvolume{3}(\bissue{22}),
\bfpage{530}
(\byear{2018}).
\doiurl{10.21105/joss.00530}
\end{barticle}
\endbibitem

\bibitem{LeClainche2017.bib}
\begin{barticle}
\bauthor{\bsnm{Le~Clainche}, \binits{S.}},
\bauthor{\bsnm{Vega}, \binits{J.M.}}:
\batitle{Higher order dynamic mode decomposition}.
\bjtitle{SIAM Journal on Applied Dynamical Systems}
\bvolume{16}(\bissue{2}),
\bfpage{882}--\blpage{925}
(\byear{2017})
\end{barticle}
\endbibitem

\bibitem{doi:10.1137/17M1125236}
\begin{barticle}
\bauthor{\bsnm{Arbabi}, \binits{H.}},
\bauthor{\bsnm{Mezić}, \binits{I.}}:
\batitle{Ergodic theory, dynamic mode decomposition, and computation of
  spectral properties of the {K}oopman operator}.
\bjtitle{SIAM Journal on Applied Dynamical Systems}
\bvolume{16}(\bissue{4}),
\bfpage{2096}--\blpage{2126}
(\byear{2017}).
\doiurl{10.1137/17M1125236}
\end{barticle}
\endbibitem

\bibitem{doi:10.1137/100813051}
\begin{barticle}
\bauthor{\bsnm{Amsallem}, \binits{D.}},
\bauthor{\bsnm{Farhat}, \binits{C.}}:
\batitle{An online method for interpolating linear parametric reduced-order
  models}.
\bjtitle{SIAM Journal on Scientific Computing}
\bvolume{33}(\bissue{5}),
\bfpage{2169}--\blpage{2198}
(\byear{2011}).
\doiurl{10.1137/100813051}
\end{barticle}
\endbibitem

\bibitem{doi:10.1137/130932715}
\begin{barticle}
\bauthor{\bsnm{Benner}, \binits{P.}},
\bauthor{\bsnm{Gugercin}, \binits{S.}},
\bauthor{\bsnm{Willcox}, \binits{K.}}:
\batitle{A survey of projection-based model reduction methods for parametric
  dynamical systems}.
\bjtitle{SIAM Review}
\bvolume{57}(\bissue{4}),
\bfpage{483}--\blpage{531}
(\byear{2015}).
\doiurl{10.1137/130932715}
\end{barticle}
\endbibitem

\bibitem{https://doi.org/10.1002/fld.2089}
\begin{barticle}
\bauthor{\bsnm{Degroote}, \binits{J.}},
\bauthor{\bsnm{Vierendeels}, \binits{J.}},
\bauthor{\bsnm{Willcox}, \binits{K.}}:
\batitle{Interpolation among reduced-order matrices to obtain parameterized
  models for design, optimization and probabilistic analysis}.
\bjtitle{International Journal for Numerical Methods in Fluids}
\bvolume{63}(\bissue{2}),
\bfpage{207}--\blpage{230}
(\byear{2010}).
\doiurl{10.1002/fld.2089}
\end{barticle}
\endbibitem

\bibitem{FarhatGrimbergManzoniQuarteroni+2020+181+244}
\begin{bchapter}
\bauthor{\bsnm{Farhat}, \binits{C.}},
\bauthor{\bsnm{Grimberg}, \binits{S.}},
\bauthor{\bsnm{Manzoni}, \binits{A.}},
\bauthor{\bsnm{Quarteroni}, \binits{A.}}:
\bctitle{Computational bottlenecks for {PROMs}: precomputation and
  hyperreduction}.
In: \bbtitle{Volume 2: Snapshot-Based Methods and Algorithms},
pp. \bfpage{181}--\blpage{244}.
\bpublisher{De Gruyter},
\blocation{Berlin, Boston}
(\byear{2021}).
\doiurl{10.1515/9783110671490-005}
\end{bchapter}
\endbibitem

\bibitem{LoiseauBruntonNoack+2020+279+320}
\begin{bchapter}
\bauthor{\bsnm{Loiseau}, \binits{J.-C.}},
\bauthor{\bsnm{Brunton}, \binits{S.L.}},
\bauthor{\bsnm{Noack}, \binits{B.R.}}:
\bctitle{From the {POD-Galerkin} method to sparse manifold models}.
In: \bbtitle{Volume 3: Applications},
pp. \bfpage{279}--\blpage{320}.
\bpublisher{De Gruyter},
\blocation{Berlin, Boston}
(\byear{2021}).
\doiurl{10.1515/9783110499001-009}
\end{bchapter}
\endbibitem

\bibitem{PanzerMohringEidLohmann+2010+475+484}
\begin{barticle}
\bauthor{\bsnm{Peuscher}, \binits{H.}},
\bauthor{\bsnm{Mohring}, \binits{J.}},
\bauthor{\bsnm{Eid}, \binits{R.}},
\bauthor{\bsnm{Lohmann}, \binits{B.}}:
\batitle{Parametric model order reduction by matrix interpolation}.
\bjtitle{Automatisierungstechnik}
\bvolume{58},
\bfpage{475}--\blpage{484}
(\byear{2010}).
\doiurl{10.1524/auto.2010.0863}
\end{barticle}
\endbibitem

\bibitem{doi:10.1137/130942462}
\begin{barticle}
\bauthor{\bsnm{Zimmermann}, \binits{R.}}:
\batitle{A locally parametrized reduced-order model for the linear frequency
  domain approach to time-accurate computational fluid dynamics}.
\bjtitle{SIAM Journal on Scientific Computing}
\bvolume{36}(\bissue{3}),
\bfpage{508}--\blpage{537}
(\byear{2014}).
\doiurl{10.1137/130942462}
\end{barticle}
\endbibitem

\bibitem{GIOVANIS2020113269}
\begin{barticle}
\bauthor{\bsnm{Giovanis}, \binits{D.G.}},
\bauthor{\bsnm{Shields}, \binits{M.D.}}:
\batitle{Data-driven surrogates for high dimensional models using gaussian
  process regression on the grassmann manifold}.
\bjtitle{Computer Methods in Applied Mechanics and Engineering}
\bvolume{370},
\bfpage{113269}
(\byear{2020}).
\doiurl{10.1016/j.cma.2020.113269}
\end{barticle}
\endbibitem

\bibitem{doi:10.2514/1.35374}
\begin{barticle}
\bauthor{\bsnm{Amsallem}, \binits{D.}},
\bauthor{\bsnm{Farhat}, \binits{C.}}:
\batitle{Interpolation method for adapting reduced-order models and application
  to aeroelasticity}.
\bjtitle{AIAA Journal}
\bvolume{46}(\bissue{7}),
\bfpage{1803}--\blpage{1813}
(\byear{2008}).
\doiurl{10.2514/1.35374}
\end{barticle}
\endbibitem

\end{thebibliography}

\end{document}